\documentclass{article}
\usepackage{epsfig,epsfig}
\usepackage{amsmath}
\usepackage{amssymb}
\usepackage{booktabs,arydshln,multirow}
\usepackage{setspace}
\usepackage{xcolor}
\usepackage{stmaryrd}
\usepackage{cite}
\newtheorem{thm}{Theorem}[section]

\newtheorem{cor}[thm]{Corollary}
\newtheorem{lemma}[thm]{Lemma}

\newtheorem{defn}[thm]{Definition}
\newtheorem{rem}[thm]{Remark}

\numberwithin{equation}{section} \topmargin=-2.5cm \oddsidemargin=0.4cm
\newcommand{\normmm}[1]{{\left\vert\kern-0.25ex\left\vert\kern-0.25ex\left\vert #1
    \right\vert\kern-0.25ex\right\vert\kern-0.25ex\right\vert}}
\begin{document}
\title{The unified theory of shifted convolution quadrature for fractional calculus
\thanks{Corresponding author.
\newline \emph{Email addresses:} mathliuyang@imu.edu.cn,
\newline \emph{Manuscript submitted to Journal~~2019}
}}
\date{ }
\author{Yang Liu$^{1*}$, Baoli Yin$^1$, Hong Li$^1$, Zhimin Zhang$^{2,3}$
\\\small{\emph{$^1$School of Mathematical Sciences, Inner Mongolia University, Hohhot 010021,
China;}}
\\\small{\emph{$^2$Beijing Computational Science Research Center, Beijing 100193, China;}}
\\\small{\emph{$^3$Department of Mathematics, Wayne State University, Detroit, MI 48202, USA}}
}
\date{}
 \maketitle
  {\color{black}\noindent\rule[0.5\baselineskip]{\textwidth}{0.5pt} }
\noindent \textbf{Abstract:} The convolution quadrature theory is a systematic approach to analyse the approximation of the Riemann-Liouville fractional operator $I^{\alpha}$ at node $x_{n}$.
In this paper, we develop the shifted convolution quadrature ($SCQ$) theory which generalizes the theory of convolution quadrature by introducing a shifted parameter $\theta$ to cover as many numerical schemes that approximate the operator $I^{\alpha}$ with an integer convergence rate as possible.
The constraint on the parameter $\theta$ is discussed in detail and the phenomenon of superconvergence for some schemes is examined from a new perspective.
For some technique purposes when analysing the stability or convergence estimates of a method applied to PDEs, we design some novel formulas with desired properties under the framework of the $SCQ$. Finally,
we conduct some numerical tests with nonsmooth solutions to further confirm our theory.
\\
 {\color{black}\noindent\rule[0.5\baselineskip]{\textwidth}{0.5pt} }
\def\REF#1{\par\hangindent\parindent\indent\llap{#1\enspace}\ignorespaces}
\newcommand{\h}{\hspace{1.cm}}
\newcommand{\hh}{\hspace{2.cm}}
\newtheorem{yl}{\hspace{0.cm}Lemma}
\newtheorem{dl}{\hspace{0.cm}Theorem}
\newtheorem{re}{\hspace{0.cm}Remark}
\renewcommand{\sec}{\section*}
\renewcommand{\l}{\langle}
\renewcommand{\r}{\rangle}
\newcommand{\be}{\begin{eqnarray}}
\newcommand{\ee}{\end{eqnarray}}
\normalsize \vskip 0.2in
\section{Introduction}
The fractional calculus has drawn much attention in recent years for its wide applications and theoretical interests, see \cite{LiChang,Jinbt1,Fordyan,Dehghan1,Fenglb1,Baleanu1,Zhaom1,LiuYBL2,LiuYBL1,Yangxj,Yangxj1,Lijch,Hesthaven}.
In this paper, we are particularly concerned about the Riemann-Liouville calculus operator $I^{\alpha}$ which is defined by
\begin{equation}\label{I.1}\begin{split}
I^{\alpha}f(x)=\frac{1}{\Gamma(\alpha)}\int_{0}^{x}(x-s)^{\alpha-1}f(s)\mathrm{d}s, \quad \text{for } \Re \alpha>0,
\end{split}
\end{equation}
and the Riemann-Liouville differential operator $I^{-\alpha} (\Re \alpha>0)$ is defined by
\begin{equation}\label{I.2}\begin{split}
I^{-\alpha}f(x)
=\frac{\mathrm{d}^n}{\mathrm{d}x^n} I^{n-\alpha}f(x)
=\frac{1}{\Gamma(n-\alpha)}\frac{\mathrm{d}^n}{\mathrm{d}x^n}\int_{0}^{x}(x-s)^{n-\alpha-1}f(s)\mathrm{d}s,
\end{split}
\end{equation}
where $n=\lceil \Re \alpha \rceil$.
When $\alpha=0$, we set $I^0=I$, the identity operator.
\par
In 1986, Lubich \cite{Lubich1} developed the convolution quadrature ($CQ$) theory to approximate the Riemann-Liouville calculus with arbitrary $\alpha \in \mathbb{C}$ at the node $x=nh$,
\begin{equation}\label{I.3}
  (I^{\alpha}f)(x) \approx h^{\alpha}\sum_{j=0}^{n}\omega_{j}f(x-jh)+h^{\alpha}\sum_{j=0}^{s}w_{n,j}f(jh),
\end{equation}
where $h$ is the mesh size of a uniform grid, $\omega_j$ denotes the convolution weight and $w_{n,j}$ is the starting weight.
For brevity, define $x_k=kh$.
For $\Re \alpha \leq 0$, the theory requires that $f^{(j)}(0)=0 (j=0,\cdots,-\lceil \Re \alpha \rceil)$.
When $\alpha$ takes integers, $I^{\alpha}$ coincides with the traditional integral operator $(\alpha>0)$ and differential operator $(\alpha<0)$, and the identity operator $(\alpha=0)$.
From this aspect, we can conclude that the $CQ$ theory generalized the traditional methods for approximating calculus operators with integer orders to fractional, i.e., arbitrary orders.
\par
Nonetheless, some classical methods are still excluded from the theory, such as the Crank-Nicolson scheme (a special case of the BDF2-$\theta$ method, see \cite{Liu1,LiuYin2}) which approximates the first derivative at the node $x_{n-\frac{1}{2}}$ and its conterpart for fractional derivative $I^{-\alpha} (\alpha \in (0,1))$ approximated at the node $x_{n-\frac{\alpha}{2}}$, with a second-order convergence rate (see the $\omega(\xi)$ of order 2 in (\ref{S.31})).
Such kind of superconvergence schemes are important for the numerical analysis for PDEs since the resulted schemes possess some good characters, see \cite{Gaogh1}.
Some other methods such as the shifted Gr$\ddot{\text{u}}$nwald formula \cite{Meerschaert}, the WSGL operator \cite{Dengwh2}, the method developed by Ding et al. \cite{Ding} and recently proposed higher-order approximation formulas \cite{Gunarathna} that generalize the fractional BDFp for $p = 1,2,\cdots,6$ are excluded from the $CQ$ theory as well, since all above methods approximate $I^{\alpha}$ at a shifted node.
\par
In this paper, we generalize the $CQ$ theory to cover the methods mentioned above by introducing a shifted parameter $\theta$ to develop the shifted convolution quadrature ($SCQ$) theory and further design some new methods with higher-order convergence rate.
Specifically, we approximate the operator $I^{\alpha} (\alpha \in \mathbb{C})$ at $x=(n-\theta)h$,
\begin{equation}\label{I.4}
  (I^{\alpha}f)(x) \approx h^{\alpha}\sum_{j=0}^{n}\omega_{j}f(x+(\theta-j)h)+h^{\alpha}\sum_{j=0}^{s}w_{n,j}f(jh).
\end{equation}
\par
The contributions of the paper are as follows
\par
$\bullet$ Develop systematic approaches to approximate the fractional calculus $I^{\alpha}f(x) (\alpha \in \mathbb{C})$ at node $x_{n-\theta}$ without assumptions on the regularity of $f(x)$.
Unlike most papers concerning the approximation at a shifted node with the parameter $\theta$, that pay little attention on the choice of $\theta$, we explore the criterion for this choice from the aspect of generating functions.
See Sec. 2.
\par
$\bullet$ Examine impacts of the parameter $\theta$ on the absolute stability regions for different numerical schemes. With a proper shifted parameter $\theta$, we can get a A-stable method, which is superior to others for some problems. See Sec. 3.
\par
$\bullet$ Construct some new approximation methods (for $I^{\alpha}f(x_{n-\theta})$) based on known approximation methods (for $I^{\alpha}f(x_{n})$), see the Theorem \ref{th.3} and Example 4.
Reveal some facts about the superconvergence for the numerical schemes (known or newly developed) from a new perspective, see the Example 2.
\par
$\bullet$ Generalize the correction technique of the $CQ$ theory to the $SCQ$ theory by introducing a parameter $\theta$, see the Theorem \ref{th.2} and Remark \ref{rem.1}.
This generalization is important since solutions of fractional calculus equations generally show some singularity at initial node.
Now with the correction technique, all methods that belong to the framework of the $SCQ$ can be modified to obtain the optimal convergence rate.
\par
$\bullet$ Apply a novel class of second-order shift-generalized Newton-Gregory formula (\ref{S.28}) to the time fractional diffusion equation with stability analysis and error estimate, see Sec. 4. Based on the analysis in foregoing sections we now pick on a set of $\theta$ with which the fractional Gr\"{o}nwall inequality can be employed in the subsequent numerical analysis.
\par
We organize the rest of the paper as follows:
In Sec. 2, we generalize the theory of the $CQ$ and develop some definitions, lemmas and theorems for the $SCQ$.
Some existing or newly proposed schemes are discussed by several examples.
In Sec. 3, we analyse the stability regions for some $SCQ$s aforementioned, discuss the impact on the regions for different $\theta$.
In Sec. 4, we devise a novel numerical scheme for the time-fractional diffusion equation by the $SCQ$ theory, with the purpose of the easy application of the discrete fractional Gr$\ddot{\text{o}}$nwall inequality.
Some lemmas are proved and the stability estimates as well as the optimal convergence order are derived. In the end of the section, we conduct some numerical tests to further confirm our theoretical analysis.
Finally, in Sec. 5, we make some conclusions and discuss some approaches we may take in the future work.
\section{Stability, consistency and convergence of the $SCQ$}
In this section, we mainly generalize the equivalence theorem developed by Lubich (Theorem 2.5, \cite{Lubich1}) which extends the classical theorem of Dahlquist \cite{Dahlquist} on linear multistep methods to fractional ones.
For convenience in the subsequent analysis we introduce some notations and definitions.
Define
\begin{equation}\label{S.1}
  (I_{h}^{\alpha}f)(x) = h^{\alpha}\sum_{j=0}^{n}\omega_{j}f(x+(\theta-j)h)+h^{\alpha}\sum_{j=0}^{s}w_{n,j}f(jh),
  \quad x=(n-\theta)h,
\end{equation}
as the shifted convolution quadrature ($SCQ$).
Denote $\Omega_{h}^{\alpha}f(x)=h^{\alpha}\sum_{j=0}^{n}\omega_{j}f(x+(\theta-j)h)$ as the convolution part and $S_{h}^{\alpha}f=h^{\alpha}\sum_{j=0}^{s}w_{n,j}f(jh)$ as the starting part of (\ref{S.1}), respectively.
Define the convolution quadrature error by
\begin{equation}\label{S.2}
 E_h^{\alpha}=\Omega_h^{\alpha}-I^{\alpha}.
\end{equation}
For the sequence of convolutions weights $\{\omega_j\}_{j=0}^{\infty}$ we associate with a generating power series $\omega(\xi)=\sum_{j=0}^{\infty}\omega_j \xi^j$, and viceversa.
\par
We note that if $f$ is continuous and $g$ is locally integrable, and $\Omega_h^{\alpha}$ is extended for $x\geq 0$ such that
\begin{equation}\label{S.3}
\Omega_{h}^{\alpha}f(x)= h^{\alpha}\sum_{0 \leq jh \leq x+\theta h}\omega_{j}f(x+(\theta-j)h),
\end{equation}
then $\Omega_h^{\alpha}$ commutes with convolution
$\Omega_h^{\alpha}(f \ast g)=(\Omega_h^{\alpha}f)\ast g$.
Hence, the convolution error $E_h^{\alpha}=\Omega_h^{\alpha}-I^{\alpha}$ satisfies (for $\Re \alpha \leq 0$ , we require $f^{(j)}(0)=0, j=0,\cdots,-\lceil \Re \alpha\rceil$)
\begin{equation}\label{S.4}
E_h^{\alpha}(f \ast g)=(E_h^{\alpha}f)\ast g.
\end{equation}
Another key property of $E_h^{\alpha}$ is the homogeneity:
\begin{equation}\label{S.5}
(E_h^{\alpha}t^{\beta-1})(x)=x^{\alpha+\beta-1}(E_{h/x}^{\alpha}t^{\beta-1})(1).
\end{equation}
We remark that (\ref{S.4}) and (\ref{S.5}) are crucial for Theorem \ref{th.1}.
The proof of (\ref{S.3}) and (\ref{S.4}) are omitted here since their correctness can be directly checked.
Next, we introduce three definitions that are closely connected in the subsequent Theorem {\ref{th.1}}:
\par
For arbitrary $\alpha \in \mathbb{C}$,
\begin{defn}\label{d.1}
A SCQ is stable for $I^{\alpha}$ if the convolution weights is
\begin{equation}\label{S.6}
\omega_n=O(n^{\alpha-1}).
\end{equation}
\end{defn}
\begin{defn}\label{d.2}
A SCQ is consistent of order $p$ for $I^{\alpha}$ if the generating function of $\{\omega_j\}_{j=0}^{\infty}$ satisfies
\begin{equation}\label{S.7}
h^{\alpha}e^{\theta h}\omega(e^{-h})=1+O(h^p).
\end{equation}
\end{defn}
\begin{defn}\label{d.3}
A SCQ is convergent of order $p$ to $I^{\alpha}$ if
\begin{equation}\label{S.8}
(E_h^{\alpha}t^{\beta-1})(1)=O(h^{\beta})+O(h^p), \quad
\text{for all } \beta \in \mathbb{C}, \beta \neq 0, -1, -2, \cdots.
\end{equation}
\end{defn}
\par
We take note of the fact that the definition of stability and convergence of the $SCQ$ are the same as the corresponding definitions of the $CQ$ (see Definition 2.1 and 2.3, \cite{Lubich1}).
The consistency of the $CQ$ is a special case of the $SCQ$ when $\theta=0$, i.e., the $CQ$ is consistent of order $p$ for $I^{\alpha}$ if $h^{\alpha}\omega(e^{-h})=1+O(h^p)$ (see Definition 2.2, \cite{Lubich1}).
\begin{rem}\label{rem.0}
Combining the homogeneity of $E_h^{\alpha}$ and the definition of convergent (\ref{S.8}), we can get
\begin{equation}\label{S.8.1}
(E_h^{\alpha}t^{\beta-1})(x_n)=O(x_n^{\alpha-1}h^{\beta})+O(x_n^{\alpha+\beta-1-p}h^p), \quad
\text{for all } \beta \in \mathbb{C}, \beta \neq 0, -1, -2, \cdots,
\end{equation}
which means for those small $\beta$ that $\alpha+\beta-1-p<0$, the convergence rate will be much lower than $p$.
Actually for fractional calculus equations, the solution is generally of weak regular at initial node.
We shall cope with such problem in the Theorem {\ref{th.2}}.
\end{rem}
\par
The following two lemmas which reveal some facts about the consistency of the $SCQ$, generalize the arguments for the $CQ$ in \cite{Lubich1}.
\begin{lemma}\label{l.1}(The counterpart of Lemma 3.1, \cite{Lubich1})
If $(E_h^{\alpha}t^{k-1})(1)=O(h^k)+O(h^p)$ for $k=1,2,3,\cdots,$ then the $SCQ$ is consistent of order $p$.
\end{lemma}
\textbf{Proof. }
First we examine the convolution error for the function $e^{t-x}$ with respect to $t$ on the interval $[0,x]$,
\begin{equation}\label{S.9}\begin{split}
e_h(x):=(E_h^{\alpha}e^{t-x})(x)&=h^{\alpha}\sum_{0 \leq jh \leq x+\theta h}\omega_{j}e^{(\theta-j)h}-(I^{\alpha}e^{t-x})(x)
\\&= h^{\alpha}e^{\theta h}\sum_{0 \leq jh \leq x+\theta h}\omega_{j}e^{-jh}-(I^{\alpha}e^{t-x})(x).
\end{split}\end{equation}
Let $x \to \infty$, the expression $h^{\alpha}e^{\theta h}\sum_{0 \leq jh \leq x+\theta h}\omega_{j}e^{-jh}$ tends to $h^{\alpha}e^{\theta h}\omega(e^{-h})$.
The rest argument of the proof is exactly the same as that of the Lemma 3.1 in \cite{Lubich1}, which is omitted here.
The proof of the lemma is completed.
\par
As is pointed out in \cite{Lubich1}, the structure of the generating function $\varpi(\xi)$ for a consistent $CQ$ is of the form (see (3.6) in \cite{Lubich1})
\begin{equation}\label{S.10}\begin{split}
\varpi(\xi)=(1-\xi)^{-\alpha}\big[c_0+c_1(1-\xi)+c_2(1-\xi)^2+\cdots+c_{N-1}(1-\xi)^{N-1}+(1-\xi)^N \tilde{r}(\xi)\big].
\end{split}\end{equation}
where $\tilde{r}(\xi)$ is holomorphic at $1$, and constants $c_j=\gamma'_j$ where $\gamma'_j$ are defined by (\ref{S.17}).
We argue that the generating function $\omega(\xi)$ for $SCQ$ can be expressed similarly by (\ref{S.10}) with a different definition of the coefficients $c_i$ that depend on $\theta$:
\begin{lemma}\label{l.2}(The counterpart of Lemma 3.2, \cite{Lubich1})
The $SCQ$ is consistent of order $p$ if and only if the coefficients $c_i$ in (\ref{S.10}) satisfy
\begin{equation}\label{S.11}\begin{split}
c_i=\gamma_i, \quad \text{for }i=0,1,\cdots,p-1,
\end{split}\end{equation}
where $\gamma_i$ are the coefficients of
\begin{equation}\label{S.12}\begin{split}
\sum_{i=0}^{\infty}\gamma_i(1-\xi)^i= \xi^{\theta}\bigg(\frac{-\ln \xi}{1-\xi}\bigg)^{-\alpha}.
\end{split}\end{equation}
\end{lemma}
\textbf{Proof.} The proof is almost the same as Lemma 3.2 in \cite{Lubich1}, and here is omitted.
\begin{rem}\label{rem.0.1}
We can construct $\omega(\xi)$ by polynomial functions.
Assume $\omega(\xi)$ has the form
\begin{equation}\label{S.2.1}
\omega(\xi)=\bigg[\frac{p_1(\xi)}{p_2(\xi)}\bigg]^{\alpha}\frac{p_3(\xi)}{p_4(\xi)},
\end{equation}
where $p_i(\xi)$ are polynomial functions, and $p_2(1)=0$.
Denote by $D$ the unit disc $|\xi|<1$, and by $\bar{D}$ the closed unit disc in the complex plane.
Then, for a stable $SCQ$, which means $\omega_n=O(n^{\alpha-1})$, it holds that $\omega(\xi)$ is analytic in $D$.
Hence, $\omega(\xi)$ can be written as
\begin{equation}\label{S.2.1.0}
\omega(\xi)=\nu(\xi)\prod_{j=0}^{\ell}(\xi-\xi_j)^{-\alpha_j},\quad \xi_0=1, \alpha_0=\alpha,
\end{equation}
where $\xi_j$ are distinct with $|\xi_j|=1$, $\nu(\xi)$ is analytic on $\bar{D}$, and $\nu(\xi_j)\neq0$, $\alpha_j\neq 0,-1,-2,\cdots$.
It can be shown that $\Re \alpha_j \leq \Re \alpha$ is equivalent to (\ref{S.6}), see \cite{Lubich1}.
We limit the choice of the parameter $\theta$ by the Condition-$\omega$,
\begin{equation}
\textbf{Condition-$\omega$: } \omega(\xi) \text{ is analytic without zeros in $D$, and } \Re \alpha_j \leq \Re \alpha.
\end{equation}
\end{rem}
\par
With the above analysis we can establish the main theorem in this paper that connects the definitions of being stable, consistent and convergent for the $SCQ$:
\begin{thm}\label{th.1}
A $SCQ$ with convolution weights defined by a generating function $\omega(\xi)$ satisfying the condition-$\omega$ is convergent of order $p$ if and only if it is stable and consistent of order $p$.
\end{thm}
\textbf{Proof.}
The theorem is a generalization of the Theorem 2.4 in \cite{Lubich1}, whose proof
consists of several lemmas.
Since the definition of stability and convergence of the $SCQ$ are the same as those of the $CQ$, we merely generalize the lemmas in \cite{Lubich1} concerning the consistency, i.e., the Lemma 3.1 and Lemma 3.2 therein.
Now with Lemma \ref{l.1} and Lemma \ref{l.2}, we complete the proof of the theorem.
\par
For a convergent $SCQ$, the following theorem shows that with the starting part $S_h^{\alpha}$, $I_h^{\alpha}f(x)$ approximates $I^{\alpha}f(x)$ uniformly for bounded $x$.
\begin{thm}\label{th.2}(See Theorem 2.4 in \cite{Lubich1})
Suppose $\omega(\xi)$ is the generating function of a convergent $SCQ$.
Then we have: For any $\beta \neq 0,-1,-2,\cdots$,
\par
(i) there exist starting weights $w_{n,j}=O(n^{\alpha-1}) (n\geq 0, j=0,\cdots,s)$ such that for any function
\begin{equation}\label{S.13}\begin{split}
f(x)=x^{\beta-1}g(x)~ \text{with } g \text{ sufficiently differentiable},
\end{split}\end{equation}
the SCQ satisfies
\begin{equation}\label{S.14}\begin{split}
I_h^{\alpha}f(x)-I^{\alpha}f(x)=O(h^p)
\end{split}\end{equation}
uniformly for $x \in [a,b]$ with $0<a<b<\infty$.
\par
(ii) there exist bounded starting weights $w_{n,j}$ such that for any function (\ref{S.13}), the estimate (\ref{S.14}) holds uniformly for bounded $x$.
\end{thm}
\begin{rem}\label{rem.1}
The starting weights $w_{n,j}$ in (i) are derived by letting
\begin{equation}\label{S.15}\begin{split}
(I_h^{\alpha}t^{q+\beta+1})(t_{n-\theta})=(I^{\alpha}t^{q+\beta+1})(t_{n-\theta})
\end{split}\end{equation}
for all integer $q \geq 0$ such that $\Re(q+\beta-1)\leq p-1$.
Note that the estimate (\ref{S.14}) requires that $x$ is bounded away from $0$.
The starting weights $w_{n,j}$ in (ii) require additionally those $\Re(q+\alpha+\beta-1)<p$ such that (\ref{S.14}) holds for bounded $x$.
\end{rem}
\par
In the rest of the section we shall collect or devise as many generating functions for the $SCQ$ as possible.
Indeed any generating function for the $CQ$ can be transformed into the one for the $SCQ$, as is proved in Theorem {\ref{th.3}}.
So let us first recall some classical or newly developed generating functions $\varpi_p(\xi)$ ($p$ is the convergence order) for the $CQ$:
\par
$\bullet$ (Theorem 2.6, \cite{Lubich1}) For an implicit linear multistep method which is stable and consistent of order $p$ with the characteristic polynomials $\rho$ and $\sigma$, if the zeros of $\sigma(\xi)$ have absolute value less than $1$, then $\varpi_p(\xi)=\big[\sigma(1/\xi)/\rho(1/\xi)\big]^{\alpha}$.
Some special cases include (see \cite{Lubich1,LiuYBL1})
\begin{equation}\label{S.16}\begin{split}
\text{The fractional trapezoidal rule: } \varpi_2(\xi)&=\frac{(1+\xi)^{\alpha}}{2^{\alpha}(1-\xi)^{\alpha}}, \quad \Re \alpha \geq 0.
\\
\text{The fractional BDF-p: } \varpi_p(\xi)&=\bigg[\sum_{i=1}^{p}\frac{1}{i}(1-\xi)^i\bigg]^{-\alpha}, \quad 1\leq p \leq 6.
\\
\text{The fractional BT-$\theta$ method: }
\varpi_2(\xi)&= \bigg[\frac{1-\theta+\theta \xi}{(3/2-\theta)-(2-2\theta)\xi+(1/2-\theta)\xi^2}\bigg]^{\alpha},
\\&\quad \text{ for }\theta \in (-\infty,1/2), \alpha<0, \text{ or }\theta \in (-\infty,1/2], \alpha\geq0.
\end{split}\end{equation}
\par
$\bullet$ The generalized Newton-Gregory formula (see \cite{Lubich1})
\begin{equation}\label{S.17}\begin{split}
\varpi_p(\xi)=(1-\xi)^{-\alpha}\big[\gamma'_0+\gamma'_1(1-\xi)+\cdots+\gamma'_{p-1}(1-\xi)^{p-1}\big],
\end{split}\end{equation}
with the coefficients $\gamma'_i$ defined by $\sum_{i=0}^{\infty}\gamma'_i(1-\xi)^i=\big(\frac{\ln \xi}{\xi-1}\big)^{-\alpha}$.
\par
$\bullet$ The fractional BN-$\theta$ method (see \cite{LiuYBL1})
\begin{equation}\label{S.18}
\begin{split}
\varpi_2(\xi)=\frac{1-\alpha\theta+\alpha\theta \xi}{\big[(3/2-\theta)-(2-2\theta)\xi+(1/2-\theta)\xi^2\big]^{\alpha}},\quad
\theta \in (-\infty,1], ~ \alpha \in \mathbb{R}, ~\alpha\theta \leq \frac{1}{2}.
\end{split}
\end{equation}

\begin{thm}\label{th.3}
Suppose $\varpi_p(\xi)$ is the generating function of a $CQ$ with convergence order $p$.
Define
\begin{equation}\label{S.19}\begin{split}
\omega(\xi)=\varpi_p(\xi)\big[\gamma''_0+\gamma''_1(1-\xi)+\cdots+\gamma''_{p-1}(1-\xi)^{p-1}\big],
\end{split}\end{equation}
where the coefficients $\gamma''_i$ satisfy
\begin{equation}\label{S.20}\begin{split}
\sum_{i=0}^{\infty}\gamma''_i(1-\xi)^i=\xi^{\theta}.
\end{split}\end{equation}
Then the SCQ with convolution weights generated by $\omega(\xi)$ is convergent of order $p$.
\end{thm}
\textbf{Proof.} For the generating function $\varpi_p(\xi)=\sum_{j=0}^{\infty}\varpi_j \xi^j$, we know that $\varpi_n=O(n^{\alpha-1})$ by the stability of the $CQ$.
Hence, we have $\omega_n=O(n^{\alpha-1})$.
By replacing $\xi$ with $e^{-h}$ in (\ref{S.20}), we get
\begin{equation}\label{S.21}\begin{split}
e^{\theta h}\sum_{i=0}^{p-1}\gamma''_i(1-e^{-h})^i=1-e^{\theta h}\sum_{i=p}^{\infty}\gamma''_i(1-e^{-h})^i=1+O(h^p).
\end{split}\end{equation}
Considering the consistency of the $CQ$, i.e.,
\begin{equation}\label{S.22}\begin{split}
h^{\alpha}\varpi_p(e^{-h})=1+O(h^p),
\end{split}\end{equation}
by combining (\ref{S.19}), (\ref{S.21}) and (\ref{S.22}), we have
\begin{equation}\label{S.23}\begin{split}
h^{\alpha}e^{\theta h}\omega(e^{-h})=1+O(h^p).
\end{split}\end{equation}
Then by Theorem \ref{th.1} the proof is completed.
\\
\par
A special case for (\ref{S.19}) that is of vital importance is when $\alpha=0$.
For any $\varpi_p(\xi)$ defined by (\ref{S.16})-(\ref{S.18}), taking $\alpha=0$, we obtain the generating function by (\ref{S.19}) for the approximation of $I^0 f(x)=f(x)$ at node $x_{n-\theta}$.
The coefficients $\gamma''_i$ in (\ref{S.20}) can be formulated as
\begin{equation}\label{S.24}\begin{split}
\gamma''_i=(-1)^i\frac{\Gamma(\theta+1)}{\Gamma(i+1)\Gamma(\theta-i+1)}.
\end{split}\end{equation}
Combining Theorem \ref{th.3}, (\ref{S.5}) and (\ref{S.8}), we have the following result
\begin{cor}\label{cor.1}
Let $f(x)=x^{\beta-1}$ with $\beta \neq 0, -1, -2, \cdots$. Then we have the following approximation formula for $f(x_{n-\theta})$ of order $p$,
\begin{equation}\label{S.25}\begin{split}
f(x_{n-\theta})-\sum_{j=0}^{p-1}\theta_{j}f(x_{n-j})=O(x_{n-\theta}^{\beta-1-p}h^p)+O(x_{n-\theta}^{-1}h^{\beta}),
\end{split}\end{equation}
where the weights $\theta_j$ satisfy
\begin{equation}\label{S.26}\begin{split}
\sum_{i=0}^{p-1}\theta_i \xi^i=\sum_{i=0}^{p-1}\gamma''_i(1-\xi)^i,
\end{split}\end{equation}
with coefficients $\gamma''_i$ defined in (\ref{S.24}).
\end{cor}
\begin{rem}\label{rem.2}
The approximation formula (\ref{S.25}) holds for $f(x)=x^{\beta-1}g(x)$ with $g(x)$ sufficiently differentiable as well.
For $p=2,3$, we can obtain two popular formulas
\begin{equation}\label{S.27}\begin{split}
f(x_{n-\theta})&\approx (1-\theta)f(x_n)+\theta f(x_{n-1}),
\\f(x_{n-\theta})&\approx  \frac{1}{2}(1-\theta)(2-\theta)f(x_n)+\theta(2-\theta)f(x_{n-1})+\frac{1}{2}\theta(\theta-1)f(x_{n-2}),
\end{split}\end{equation}
see also (17), (18) in \cite{Dimitrov}.
Considering the condition-$\omega$, we limit $\theta$ for (\ref{S.27}) to satisfying $\theta \leq \frac{1}{2}$ and $\theta \leq 1-\frac{\sqrt{2}}{2}$, respectively.
As one can see from (\ref{S.25}), to get a convergence order of $p$ for bounded $x$ if $f(x)$ is not so regular, the staring part is needed according to the Theorem \ref{th.2}.
\end{rem}
\par
Another family of generating functions for the $SCQ$ is the shift-generalized Newton-Gregory formula which is a natural result from the analysis of Lemma \ref{l.2}:
\begin{cor}\label{cor.2}(The shift-generalized Newton-Gregory formula)
The SCQ with weights generated by the following generating function is convergent (to $I^{\alpha}$) of order $p$,
\begin{equation}\label{S.28}\begin{split}
\omega(\xi)=(1-\xi)^{-\alpha}\big[\gamma_0+\gamma_1(1-\xi)+\cdots+\gamma_{p-1}(1-\xi)^{p-1}\big],
\end{split}\end{equation}
where $\gamma_i$ are the coefficients of
\begin{equation}\label{S.29}\begin{split}
\sum_{i=0}^{\infty}\gamma_i(1-\xi)^i= \xi^{\theta}\bigg(\frac{-\ln \xi}{1-\xi}\bigg)^{-\alpha}.
\end{split}\end{equation}
\end{cor}
One can find out that the shift-generalized Newton-Gregory formula reduces to (\ref{S.17}) when $\theta=0$.
We conclude this section by further exploring some numerical methods that approximate $f(x)$ at node $x_{n-\theta}$ in the following examples. All methods mentioned in the examples belong to the SCQ, hence correction technique can be applied to the methods if the solution is not regular enough.
For simplicity, we assume $\alpha \in \mathbb{R}$.
\par
\textbf{Example 1. }(Generalized shifted Gr$\ddot{\text{u}}$nwald formula)
\par
Expanding (\ref{S.29}), we can easily derive that
\begin{equation}\label{S.30}\begin{split}
\gamma_0=1,\quad
\gamma_1=-\frac{\alpha}{2}-\theta, \quad \gamma_2=\frac{1}{8}(\alpha+2\theta)^2-\frac{1}{24}(5\alpha+12\theta).
\end{split}\end{equation}
For the case with convergence order $p=1$, we have the approximation formula as
\begin{equation}\label{S.30.1}\begin{split}
(I^{\alpha}f)(x) \approx h^{\alpha}\sum_{j=0}^{n}\omega_{j}f\big(x+(\theta-j)h\big)+h^{\alpha}\sum_{j=0}^{s}w_{n,j}f(jh), \quad
\text{with } \omega(\xi)=(1-\xi)^{-\alpha}.
\end{split}\end{equation}
Now by assuming that $f$ is sufficiently smooth (hence, the starting part can be omitted), and that $\theta$ takes nonnegative integers, the relation  (\ref{S.30.1}) reduces to the shifted Gr$\ddot{\text{u}}$nwald formula (see \cite{Meerschaert} with $\alpha$ replaced by $-\alpha$).
Recently, Chen et al. \cite{ChenMartin} has applied the Gr$\ddot{\text{u}}$nwald formula ($\theta=0$) to the time fractional PDEs and derived a sharp convergence rate which is in line with (\ref{S.8.1}).
\par
\textbf{Example 2. }(Discussion on the superconvergence)
\par
With the coefficients in (\ref{S.30}), the generating function $\omega(\xi)$ for the shift-generalized Newton-Gregory formula of order $p=2,3$ can be formulated as
\begin{equation}\label{S.31}\begin{split}
\omega(\xi)=(1-\xi)^{-\alpha}\bigg[1-\bigg(\frac{\alpha}{2}+\theta\bigg)&(1-\xi)\bigg]
 \quad \text{of order } 2, \text{ with }\theta \leq \frac{1-\alpha}{2},
\\
\omega(\xi)=(1-\xi)^{-\alpha}\bigg\{1-\bigg(\frac{\alpha}{2}+\theta\bigg)(1-\xi)+\bigg[\frac{1}{8}&(\alpha+2\theta)^2-\frac{1}{24}(5\alpha+12\theta)\bigg](1-\xi)^2\bigg\} \quad \text{of order } 3.
\end{split}\end{equation}
Both of the formula were also proposed by Dimitrov \cite{Dimitrov} by using the theory developed in \cite{Tadjeran}.
An obvious conclusion is that $\omega(\xi)$ defined by (\ref{S.31}) which is of order $2$ shows some superconvergence when $\theta=-\frac{\alpha}{2}$ (see also \cite{Gaogh1}, with $\alpha$ replaced by $-\alpha$), since generally a $CQ$ with the generating function $(1-\xi)^{-\alpha}$ approximates $I^{\alpha}f(x_n)$ with a lower order.
From the perspective of generating function (\ref{S.28}), we can always find some superconvergence points that effective numerical methods can be proposed.
To the best knowledge of the authors, there is no literature exploring the superconvergence property of the scheme with $\omega(\xi)$ defined by (\ref{S.31}) which is of order $3$, that from the discussion above, a $SCQ$ with $\omega(\xi)$ defined by
\begin{equation}\label{S.31.1}\begin{split}
\omega(\xi)=(1-\xi)^{-\alpha}\bigg[1-\bigg(\frac{\alpha}{2}+\theta\bigg)&(1-\xi)\bigg]
\end{split}\end{equation}
is convergent of order $3$ provided that
\begin{equation}\label{S.31.2}\begin{split}
\theta=\frac{1-\alpha}{2}-\frac{1}{2}\sqrt{1-\frac{\alpha}{3}}, \quad
\text{for } \alpha \leq 3.
\end{split}\end{equation}
One can check that this choice of $\theta$ satisfies the condition-$\omega$.
\par
\textbf{Example 3. }(The WSGL method, see \cite{Dengwh2,Zeng1,WangZB,Liuydyw2,Duywly,LiuZhangmin,SunCC,ChenDeng})
\par
Assume two integers $p > q$, and let $\theta=p$, consider the following generating function
\begin{equation}\label{S.32}\begin{split}
\omega(\xi)=(1-\xi)^{-\alpha}\bigg[\frac{-\alpha-2q}{2(p-q)}+\frac{2p+\alpha}{2(p-q)}\xi^{p-q}\bigg],
\end{split}\end{equation}
whose coefficients can be formulated as
\begin{equation}\label{S.33}\begin{split}
\omega_k=\frac{-\alpha-2q}{2(p-q)}\varrho_k+\frac{2p+\alpha}{2(p-q)}\varrho_{k-(p-q)},
\end{split}\end{equation}
where $\varrho_k$ are the coefficients of $(1-\xi)^{-\alpha}$ and define $\varrho_k=0$ for $k<0$.
Considering the condition-$\omega$, we assume $(p,q)$ satisfies
\begin{equation}\label{S.33.1}
p+q+\alpha \leq 0.
\end{equation}
One can easily check that the $SCQ$ with $\omega(\xi)$ defined in (\ref{S.32}) is the WSGL operator (see \cite{Dengwh2}, with $\alpha$ replaced by $-\alpha$), which is convergent of order $2$ by Theorem \ref{th.1}.
Actually, by (\ref{S.33}), we have $\omega_n=O(n^{\alpha-1})$ since $\varrho_n=O(n^{\alpha-1})$ (see the proof of Lemma \ref{l.4} with $\alpha$ replaced by $-\alpha$).
And, by the Taylor expansion formulas, we can get
\begin{equation}\label{S.33.2}\begin{split}
h^{\alpha}e^{p h}\omega(e^{-h})&=h^{\alpha}e^{p h}(1-e^{-h})^{-\alpha}\bigg[\frac{-\alpha-2q}{2(p-q)}+\frac{2p+\alpha}{2(p-q)}e^{h(q-p)}\bigg]
\\&=
h^{\alpha}(1-e^{-h})^{-\alpha}\bigg[\frac{-\alpha-2q}{2(p-q)}e^{hp}+\frac{2p+\alpha}{2(p-q)}e^{hq}\bigg]
\\&=
\bigg[1-\frac{h}{2}+O(h^2)\bigg]^{-\alpha}\bigg[1-\frac{\alpha}{2}h+O(h^2)\bigg]
\\&=
\bigg[1+\frac{\alpha}{2}h+O(h^2)\bigg]\bigg[1-\frac{\alpha}{2}h+O(h^2)\bigg]
\\&=1+O(h^2),
\end{split}\end{equation}
which implies the WSGL operator is stable and consistent.
Actually, the second-order WSGL method is constructed by the first-order fractional BDF with specially designed weights. For some other numerical methods constructed by the fractional BDF but with higher-order convergence rates, see \cite{Dengwh2,LiuZhangmin,SunCC,ChenDeng}.
\par
\textbf{Example 4. }(Further discussion on the WSGL method)
\par
A interesting consideration is that with the structure of the generating function (\ref{S.32}), i.e., with
\begin{equation}\label{S.34}\begin{split}
\omega(\xi)=(1-\xi)^{-\alpha}(\lambda_1+\lambda_2 \xi^{m}),
\quad  \lambda_1 \geq |\lambda_2|,
\end{split}\end{equation}
where $m$ is a positive integer, can we propose a new formula that is convergent (to $I^{\alpha}$) of order $3$?
Actually, with the following coefficients
\begin{equation}\label{S.35}\begin{split}
\lambda_1=1-\lambda_2, \quad
\lambda_2=\frac{\alpha+2\theta}{2m},\quad
m=1+\frac{\alpha+2\theta}{2}-\frac{5\alpha+12\theta}{6(\alpha+2\theta)},
\end{split}\end{equation}
one can easily check that the $SCQ$ with $\omega(\xi)$ defined by (\ref{S.35}) is stable and consistent of order $3$.
Furthermore, since $\alpha+2\theta$ and $m$ cannot be zero by (\ref{S.35}), it seems that there is no superconvergence point for the $(1-\xi)^{-\alpha}$ to devise a formula of convergence order $3$.
For the application of (\ref{S.34}), by considering the condition $\lambda_1 \geq |\lambda_2|$ (which is derived by the condition-$\omega$), we can take
\begin{equation}\label{S.36}\begin{split}
\theta=\frac{m-\alpha}{2}-\frac{1}{2}\sqrt{m^2-\frac{\alpha}{3}},\quad
\text{for } \alpha \leq 3m^2.
\end{split}\end{equation}
We remark here that (\ref{S.34}) generalizes (\ref{S.31.1}) from some aspects that, if we take $m=1$, under the condition $(\ref{S.36})$ (which reduces to (\ref{S.31.2})), then the generating function $\omega(\xi)$ defined by (\ref{S.34}) with (\ref{S.35}) is exactly the same as the $\omega(\xi)$ defined by (\ref{S.31.1}).
\par
\textbf{Example 5. }(Generalized BDF2-$\theta$ method)
\par
In this example we consider the generalized BDF2-$\theta$ method that generalizes the work of Liu et al. \cite{Liu1} and Ding et al. \cite{Ding}.
Define the generating function $\omega(\xi)$ by
\begin{equation}\label{S.37}\begin{split}
\omega(\xi)=\bigg(\frac{3\alpha+2\theta}{2\alpha}-\frac{2\alpha+2\theta}{\alpha}\xi+\frac{\alpha+2\theta}{2\alpha}\xi^2\bigg)^{-\alpha},
\quad \text{with }\alpha(\alpha+\theta)\geq 0,
\end{split}\end{equation}
which can be reformulated as
\begin{equation}\label{S.38}\begin{split}
\omega(\xi)=\bigg(\frac{3}{2}+\frac{\theta}{\alpha}\bigg)^{-\alpha}(1-\xi)^{-\alpha}\bigg(1-\frac{\alpha+2\theta}{3\alpha+2\theta}\xi\bigg)^{-\alpha}.
\end{split}\end{equation}
The generating function (\ref{S.37}) has been proposed by Gunarathna et al. \cite{Gunarathna} for the fractional derivative, i.e., for $\alpha<0$.
For the case $\theta$ takes integers and other higher-order formulas, see \cite{LiCai}.
By careful derivation, we can check that the $SCQ$ with $\omega(\xi)$ defined in (\ref{S.37})  is convergent of order $2$.
If we take $\alpha=-1$, (\ref{S.37}) is reduced to the BDF-$\theta$ method (see \cite{Liu1}) with the generating function as,
\begin{equation}\label{S.39}\begin{split}
\omega(\xi)=\frac{3-2\theta}{2}-(2-2\theta)\xi+\frac{1-2\theta}{2}\xi^2.
\end{split}\end{equation}
If we take $\theta=-\frac{1}{2}$, we get a $SCQ$ that approximates $f$ at the node $x_{n+\frac{1}{2}}$ (see \cite{Ding} with $\alpha$ replaced by $-\alpha$),
\begin{equation}\label{S.40}\begin{split}
\omega(\xi)=\bigg(\frac{3\alpha-1}{2\alpha}-\frac{2\alpha-1}{\alpha}\xi+\frac{\alpha-1}{2\alpha}\xi^2\bigg)^{-\alpha}.
\end{split}\end{equation}
There are two points deserve discussion:
\\
a) The generating function (\ref{S.37}) also implies that some superconvergence properties at node $x_{n+\frac{\alpha}{2}}$ by taking $\theta=-\frac{\alpha}{2}$.
\\
b) By taking $\theta=-\alpha$, we can obtain a shorter or simpler generating function
\begin{equation}\label{S.41}\begin{split}
\omega(\xi)=2^{\alpha}(1-\xi^2)^{-\alpha},
\end{split}\end{equation}
and we call the corresponding $SCQ$ the fractional central difference method, for the reason that if we take $\alpha=-1$, then the generating function $\omega(\xi)=\frac{1}{2}(1-\xi^2)$ implies the classical central difference scheme
\begin{equation}\label{S.42}\begin{split}
f'(x_{n-1}) \approx \frac{f(x_n)-f(x_{n-2})}{2}.
\end{split}\end{equation}
However, the application of the fractional central difference method is limited.
See table \ref{tab4} in section 3 and the notation therein.
\section{Stability regions}
In this section we pay special attention on the stability regions (see Definition {\ref{d.5}}) of the $SCQ$ when applied to the fractional equations after omitting the starting part.
This work is motivated by the fact that the condition-$\omega$ can not guarantee a $SCQ$ is A($\frac{\pi}{2}$)-stable (see Definition {\ref{d.5}}).
As is well known, for some problems numerical schemes of A($\frac{\pi}{2}$)-stable are superior to the others.
To be specific, we analyse the following two models
\begin{equation}\label{R.1}\begin{split}
{}_C D_{0,x}^{\alpha} y(x)=\lambda y(x), \quad y(0)=y_0, \quad \alpha\in(0,1),
\end{split}\end{equation}
and
\begin{equation}\label{R.2}\begin{split}
I^{\alpha} y(x)=\frac{1}{\lambda}y(x),\quad \alpha\in(0,1),
\end{split}\end{equation}
respectively.
The operator ${}_C D_{0,x}^{\alpha}$ denotes the Caputo fractional derivative operator of order $\alpha$, which is defined by
\begin{equation}\label{R.3}\begin{split}
{}_C D_{0,x}^{\alpha} y(x)=\frac{1}{\Gamma(1-\alpha)}\int_{0}^{x}\frac{y'(s)}{(x-s)^{\alpha}}\mathrm{d}s,\quad
\alpha \in (0,1).
\end{split}\end{equation}
With the relation
\begin{equation}\label{R.4}\begin{split}
{}_C D_{0,x}^{\alpha} y(x)=I^{-\alpha}(y(x)-y_0), \quad \alpha \in (0,1),
\end{split}\end{equation}
we can approximate the Caputo type derivatives by the $SCQ$ developed in the Sec. 2.
\par
Define $z:=\lambda h^{\alpha}$, the numerical schemes for equation (\ref{R.1}) and (\ref{R.2}) are as follows
\begin{equation}\label{R.5}\begin{split}
\text{Scheme-I: }&\sum_{k=0}^{n}\omega_{n-k}y^k=z\sum_{j=0}^{p-1}\theta_j y^{n-j}+f(x_{n-\theta}),
\\
\text{Scheme-II: }&z\sum_{k=0}^{n}\omega_{n-k}y^k=\sum_{j=0}^{p-1}\theta_j y^{n-j},
\end{split}\end{equation}
respectively, where $\theta_j$ are defined by (\ref{S.26}), and
\begin{equation}\label{R.6}\begin{split}
f(x_{n-\theta}):=h^{\alpha}(I^{-\alpha}y_0)(x_{n-\theta})=\frac{y_0}{(n-\theta)^{\alpha}\Gamma(1-\alpha)}.
\end{split}\end{equation}
\par
Before examine the stability region of the numerical schemes (\ref{R.5}), we state some properties of the analytic solutions of the model equation (\ref{R.1}) and (\ref{R.2}).
\begin{defn}\label{d.4} (See \cite{Podlubny})
The two-parameter Mittag-Leffler function $E_{\alpha,\beta}(x)$ is defined by,
\begin{equation}\label{R.7}\begin{split}
E_{\alpha,\beta}(x):=\sum_{j=0}^{\infty}\frac{x^j}{\Gamma(j\alpha+\beta)},
\text{ for those $x$ that make the series converge}.
\end{split}\end{equation}
\end{defn}
For simplicity, we define $E_{\alpha}(x):=E_{\alpha,1}(x)$.
\begin{lemma}\label{l.3}
For $|\arg \lambda -\pi|<(1-\frac{\alpha}{2})\pi$, the solutions of (\ref{R.1}) and (\ref{R.2}) satisfy
\begin{equation}\label{R.8}\begin{split}
y(x)\to 0 \quad \text{as}\quad x \to \infty.
\end{split}\end{equation}
\end{lemma}
\textbf{Proof. }
For the equation (\ref{R.2}), which is a special case of (1.2) in \cite{Lubich2}, see the proof therein.
We mainly focus on the differential equation (\ref{R.1}).
By Laplace transform, we can express the analytic solution by the Mittag-Leffler function $E_{\alpha}(x)$,
\begin{equation}\label{R.9}\begin{split}
y(x)=y_0E_{\alpha}(\lambda x^{\alpha}).
\end{split}\end{equation}
Note that for $\alpha \in (0,2)$, we have the asymptotic property (see Theorem 1.6, p.35, \cite{Podlubny}) for $E_{\alpha}(x)$ that with $\frac{\alpha}{2}\pi < \mu < \min\{\pi,\alpha \pi\}$,
\begin{equation}\label{R.10}\begin{split}
|E_{\alpha}(x)|\leq \frac{C}{1+|x|}, \quad \mu \leq |\arg x|\leq \pi, \quad |x|\geq 0,
\end{split}\end{equation}
where $C$ is a real constant.
By replacing $x$ of $E_{\alpha}(x)$ with $\lambda x^{\alpha}$, the proof of the lemma is completed.
\par
The Lemma \ref{l.3} naturally leads to the following definition (see also \cite{Lubich2}),
\begin{defn}\label{d.5}
The stability region $S$ of a $SCQ$ is the set of all complex $z=\lambda h^{\alpha}$ for which the numerical solutions $y^n$ of (\ref{R.5}) satisfy
\begin{equation}\label{R.11}\begin{split}
y^n \to 0 \quad \text{as} \quad n \to \infty.
\end{split}\end{equation}
Further, we call a numerical method $A(\delta)$-stable if $S$ contains the sector $|\arg z-\pi|<\delta$.
\end{defn}
\begin{thm}\label{th.4}
The stability regions of the numerical schemes (\ref{R.5}) are
\begin{equation}\label{R.12}\begin{split}
\text{For Scheme-I: }\quad&\mathbb{C}\setminus\bigg\{z:z=\omega(\xi)\bigg/\sum_{j=0}^{p-1}\theta_j\xi^j, |\xi| \leq 1\bigg\},
\\
\text{For Scheme-II: }\quad&
\mathbb{C}\setminus\bigg\{z:z=\sum_{j=0}^{p-1}\theta_j\xi^j\bigg/\omega(\xi), |\xi| \leq 1\bigg\},
\end{split}\end{equation}
respectively.
\end{thm}
\text{Proof. }The technique used in this theorem is the same as the Theorem 2.1 in \cite{Lubich2}. We omit the proof here.
\par
We conclude this section by illustrating the stability regions of some $SCQ$s and make some notations.
Generally, we require that there exists a positive $x_0$ such that the interval $(-x_0,0)$ is contained in the stability region $S$.
\par
Table \ref{tab1} illustrates the stability regions for the shift-generalized Newton-Gregory formula of order $2$. For the Scheme-I with $\theta \leq \frac{1}{2}$, the method is $A(\frac{\pi}{2})$-stable.
If $\theta$ satisfies $\frac{1}{2}<\theta<\frac{\alpha+1}{2}$, the method is conditional stable which means for arbitrary $\lambda < 0$, the step size $h$ must be small enough ($\lambda h^{\alpha}>-x_0$) such that for fixed $h$, the solution $y^n$ tends to $0$ as $n$ tends to infinity.
If $\theta=\frac{\alpha+1}{2}$, then for any $x_0>0$, $(-x_0,0) \cap S= \varnothing$, which means for any step size $h>0$, $y^n$ oscillates or blows up as $n$ tends to infinity.
For the scheme-II, by condition-$\omega$, we require $\theta \leq \frac{1-\alpha}{2}$, in which the method is $A(\frac{\pi}{2})$-stable.
If $\theta$ exceeds this value, we will get a conditional stable method, provided $\theta < \frac{1}{2}$.
\par
In Table \ref{tab2}, we briefly depict the stability regions for WSGL operators, with the choice of pairs $(p,q)=(0,-2)$ or $(1,-2)$ for the Scheme-I, and $(p,q)=(0,-2)$ or $(-1,-2)$ for the Scheme-II.
One can find that only with $(p,q)=(1,-2)$ the method is conditional stable and for other three pairs the method is $A(\frac{\pi}{2})$-stable.
\par
In Table \ref{tab3}, we show some stability regions for the $SCQ$ with $\omega(\xi)$ defined in (\ref{S.37}).
For the Scheme-I under the condition-$\omega$, we can conclude that with $\theta \leq \min\{\alpha, \frac{1}{2}\}$, the method is $A(\frac{\pi}{2})$-stable.
For the Scheme-II, the method is $A(\frac{\pi}{2})$-stable if $-\alpha \leq \theta \leq \frac{1}{2}$.
An interesting phenomenon is that when $\theta>\frac{1}{2}$, there exists a $x_0>0$, such that for any $h>x_0$, the solution $y^n\to 0 (n\to \infty)$.
\par
\begin{figure}[h]
\begin{center}
\begin{minipage}{15.5cm}
  \centering\includegraphics[width=15cm]{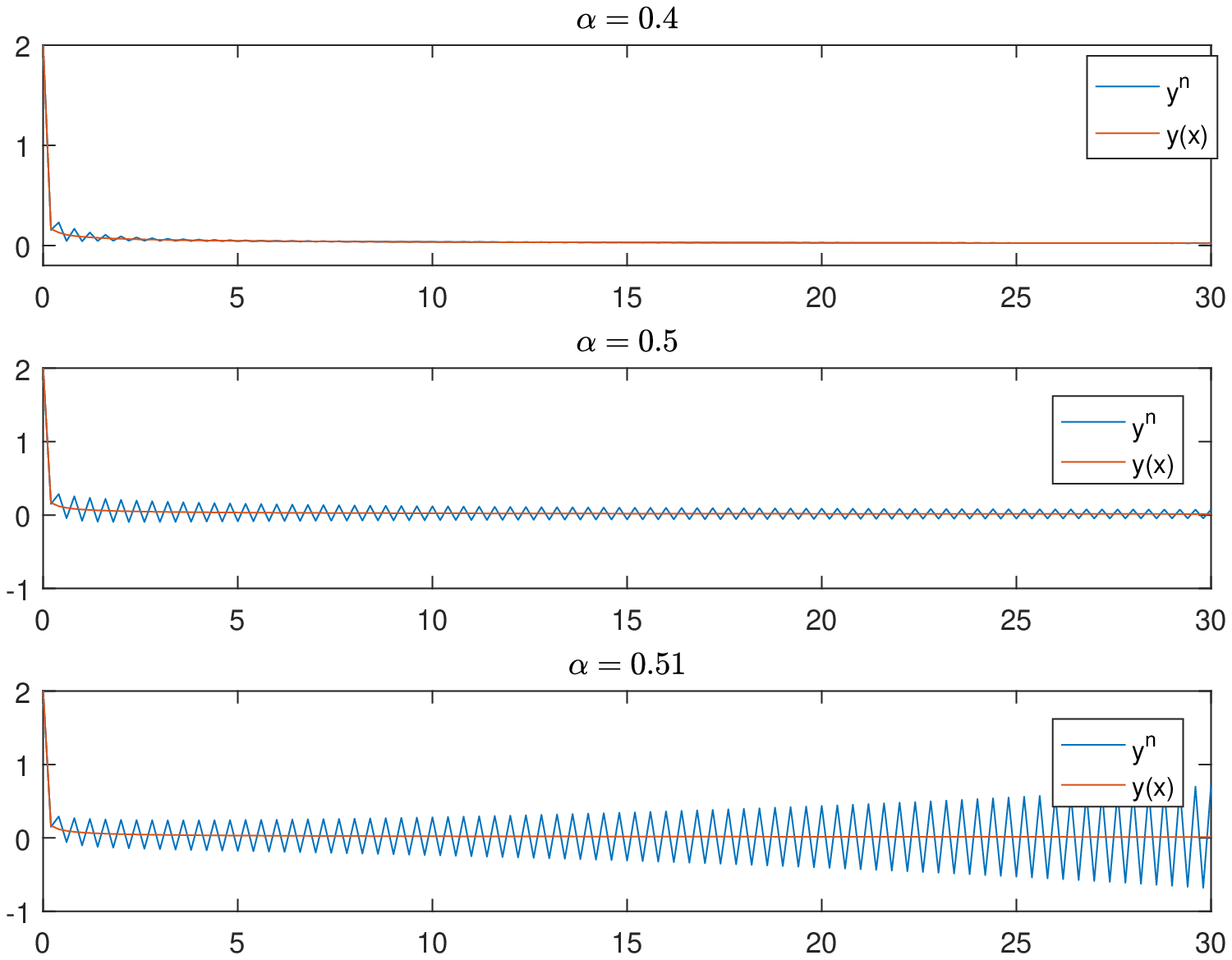}
  \caption{$h=0.2$, $\lambda=-15$.}\label{C15}
\end{minipage}
\end{center}
\end{figure}

In Table \ref{tab4} we analyse the reason why the application for the fractional central difference method is limited.
For the Scheme-I, one can find that the solution $y^n$ blows up for $\alpha > \frac{1}{2}$, since for any $x_0>0$, $(-x_0, 0) \cap S = \varnothing$.
See also Fig. \ref{C15} for the numerical solutions of the Scheme-I with $\lambda=-15$, $h=0.2$.
Nonetheless, the method is quite suitable for the Scheme-II which is $A(\frac{\pi}{2})$-stable for any $\alpha \in (0,1)$.
\begin{table}[h]
\centering
 \caption{The shift-generalized Newton-Gregory formula of order $2$}\label{tab1}
{\renewcommand{\arraystretch}{2}
\renewcommand{\tabcolsep}{0cm}
\begin{tabular}{|c|c|}
\hline
Scheme-I & Scheme-II \\
\hline
$\omega(\xi)=(1-\xi)^{\alpha}\big[1+\big(\frac{\alpha}{2}-\theta\big)(1-\xi)\big]$ & $\omega(\xi)=(1-\xi)^{-\alpha}\big[1-\big(\frac{\alpha}{2}+\theta\big)(1-\xi)\big]$ \\
\hline
\begin{minipage}{7.5cm}
  \centering\includegraphics[width=7cm]{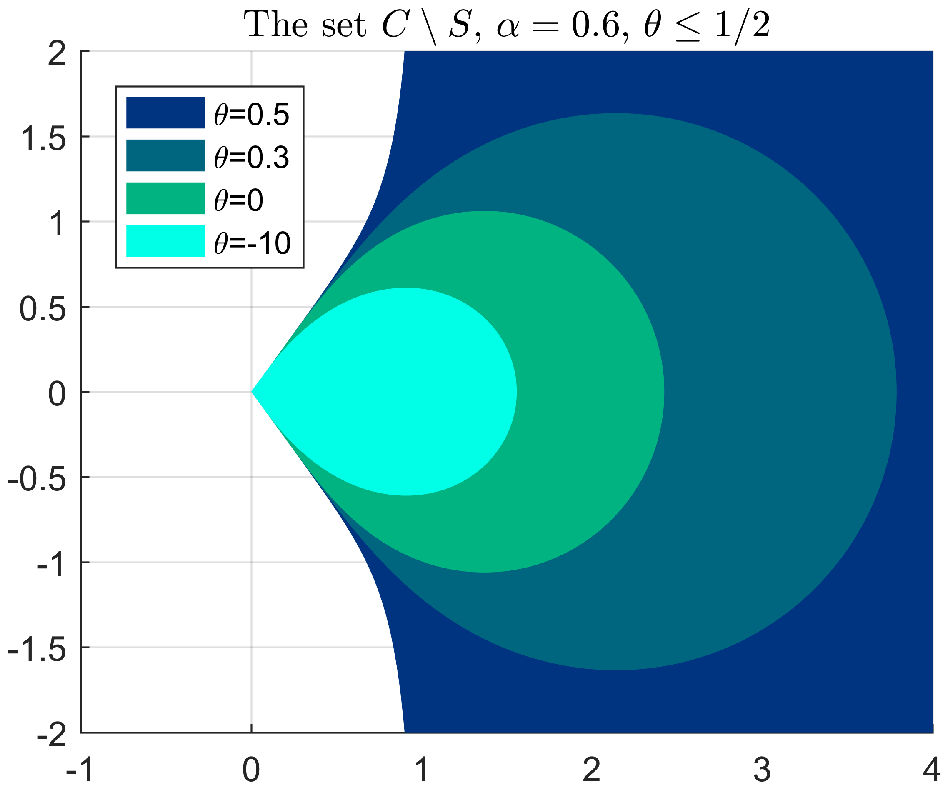}\label{C1}
\end{minipage}
&
\begin{minipage}{7.5cm}
  \centering\includegraphics[width=7cm]{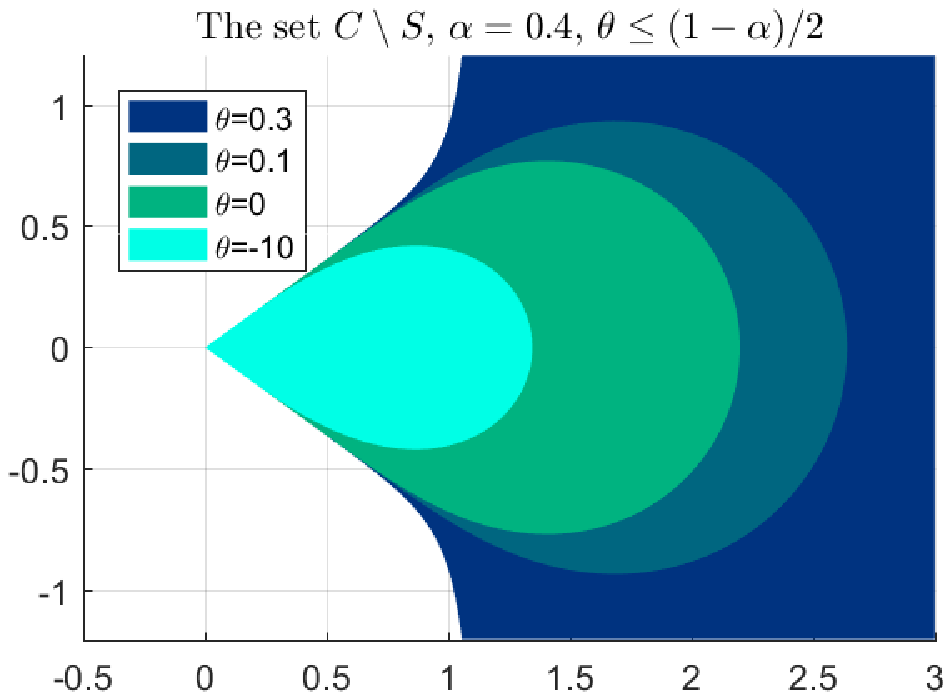}\label{C2}
\end{minipage} \\
\begin{minipage}{7.5cm}
  \centering\includegraphics[width=7cm]{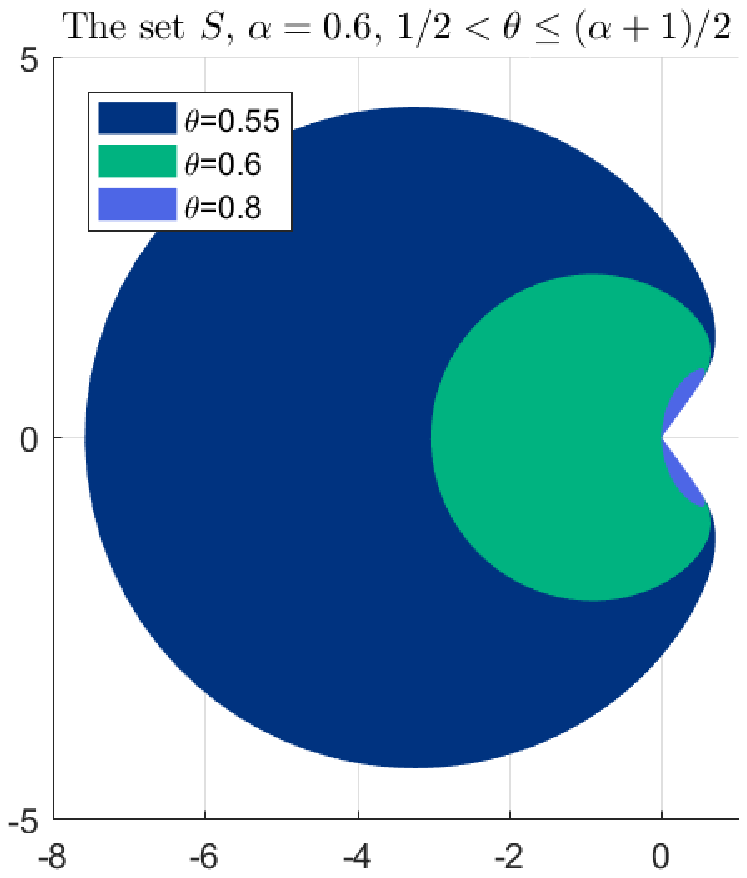}\label{C3}
\end{minipage}
&
\begin{minipage}{7.5cm}
  \centering\includegraphics[width=7cm]{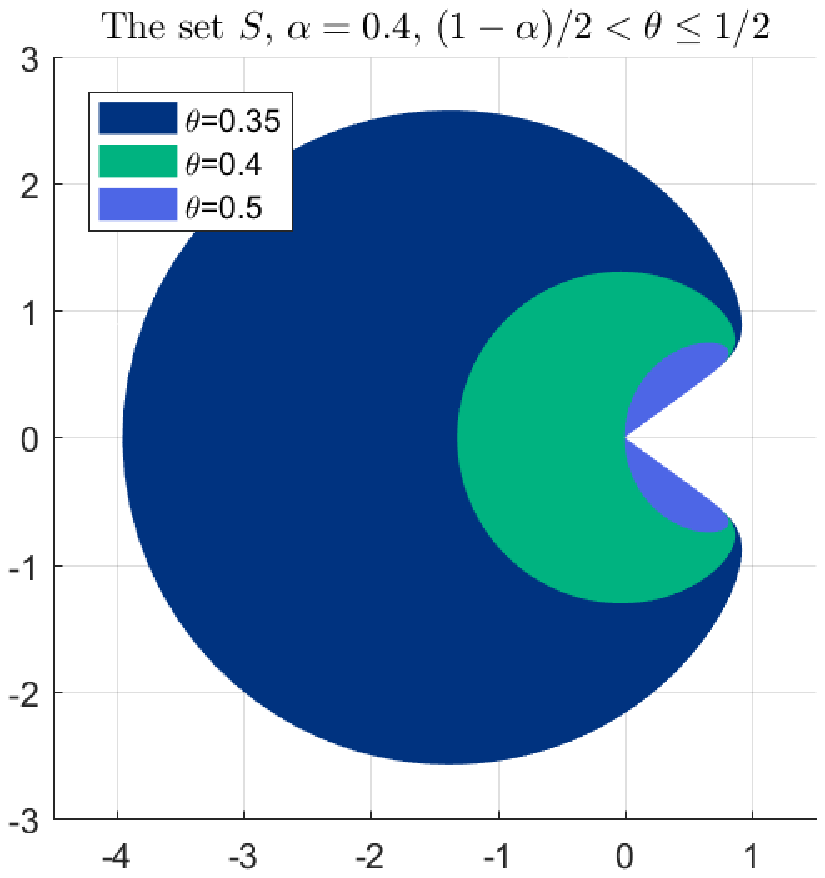}\label{C4}
\end{minipage}\\
\hline
\end{tabular}}
\end{table}
\begin{table}[h]
\centering
 \caption{The WSGL operator with $q=-2, \theta=p$}\label{tab2}
{\renewcommand{\arraystretch}{2}
\renewcommand{\tabcolsep}{0cm}
\begin{tabular}{|c|c|}
\hline
Scheme-I & Scheme-II \\
\hline
$\omega(\xi)=(1-\xi)^{\alpha}\big[\frac{4+\alpha}{2(p+2)}+\frac{2p-\alpha}{2(p+2)}\xi^{p+2}\big]$ & $\omega(\xi)=(1-\xi)^{-\alpha}\big[\frac{4-\alpha}{2(p+2)}+\frac{2p+\alpha}{2(p+2)}\xi^{p+2}\big]$ \\
\hline
\begin{minipage}{7.5cm}
  \centering\includegraphics[width=7cm]{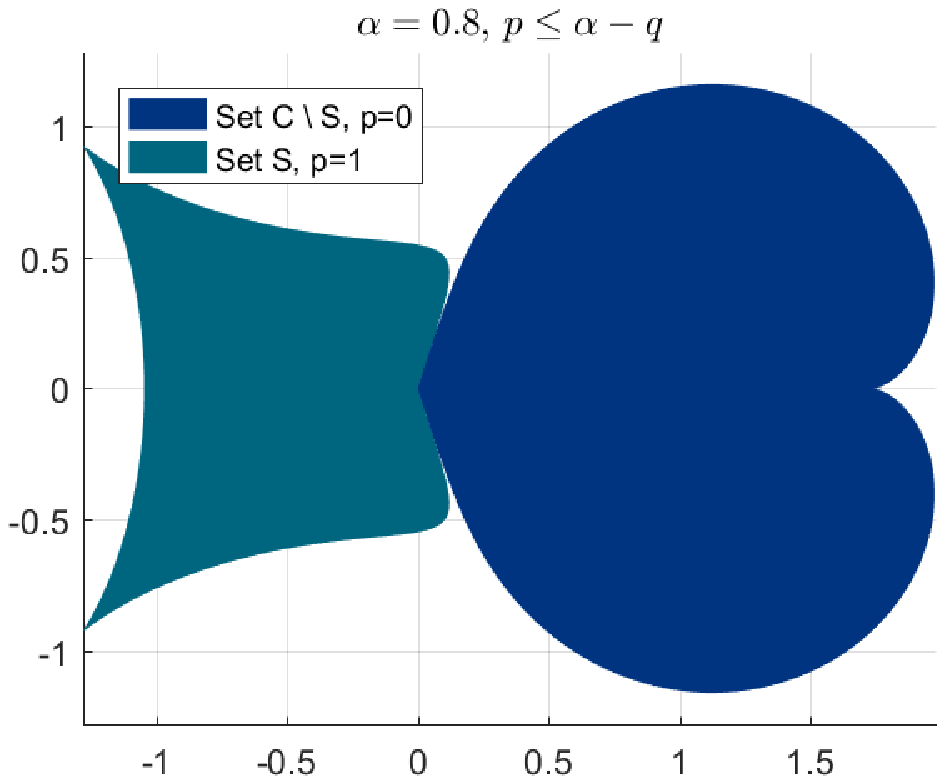}\label{C5}
\end{minipage}
 &
\begin{minipage}{7.5cm}
  \centering\includegraphics[width=7cm]{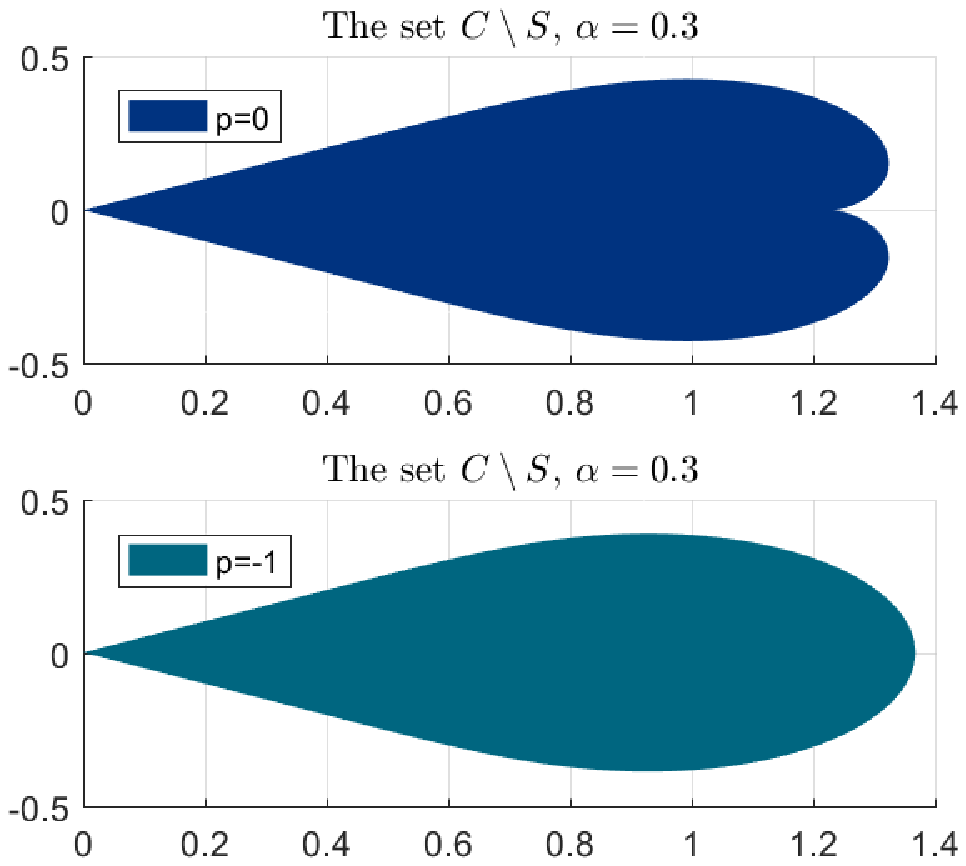}\label{C6}
\end{minipage}
\\
\hline
\end{tabular}}
\end{table}
\begin{table}[h]
\centering
 \caption{The $SCQ$ with $\omega(\xi)$ in (\ref{S.37})}\label{tab3}
{\renewcommand{\arraystretch}{2}
\renewcommand{\tabcolsep}{0cm}
\begin{tabular}{|c|c|}
\hline
Scheme-I & Scheme-II \\
\hline
$\omega(\xi)=\big(\frac{3\alpha-2\theta}{2\alpha}-\frac{2\alpha-2\theta}{\alpha}\xi+\frac{\alpha-2\theta}{2\alpha}\xi^2\big)^{\alpha}$ & $\omega(\xi)=\big(\frac{3\alpha+2\theta}{2\alpha}-\frac{2\alpha+2\theta}{\alpha}\xi+\frac{\alpha+2\theta}{2\alpha}\xi^2\big)^{-\alpha}$ \\
\hline
\begin{minipage}{7.5cm}
  \centering\includegraphics[width=7cm]{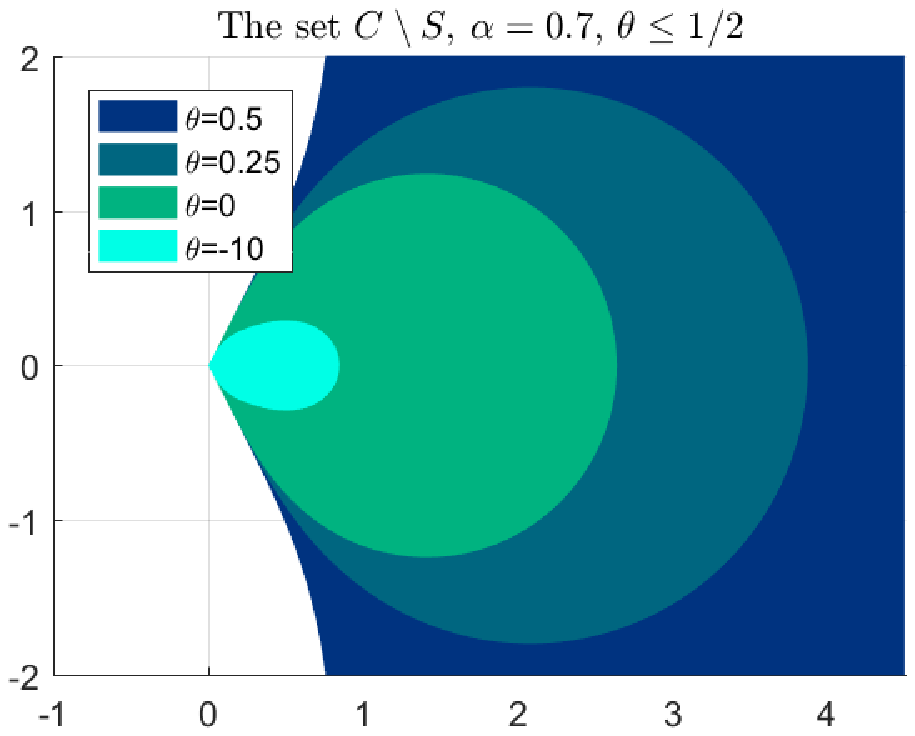}\label{C7}
\end{minipage}
&\begin{minipage}{7.5cm}
  \centering\includegraphics[width=7cm]{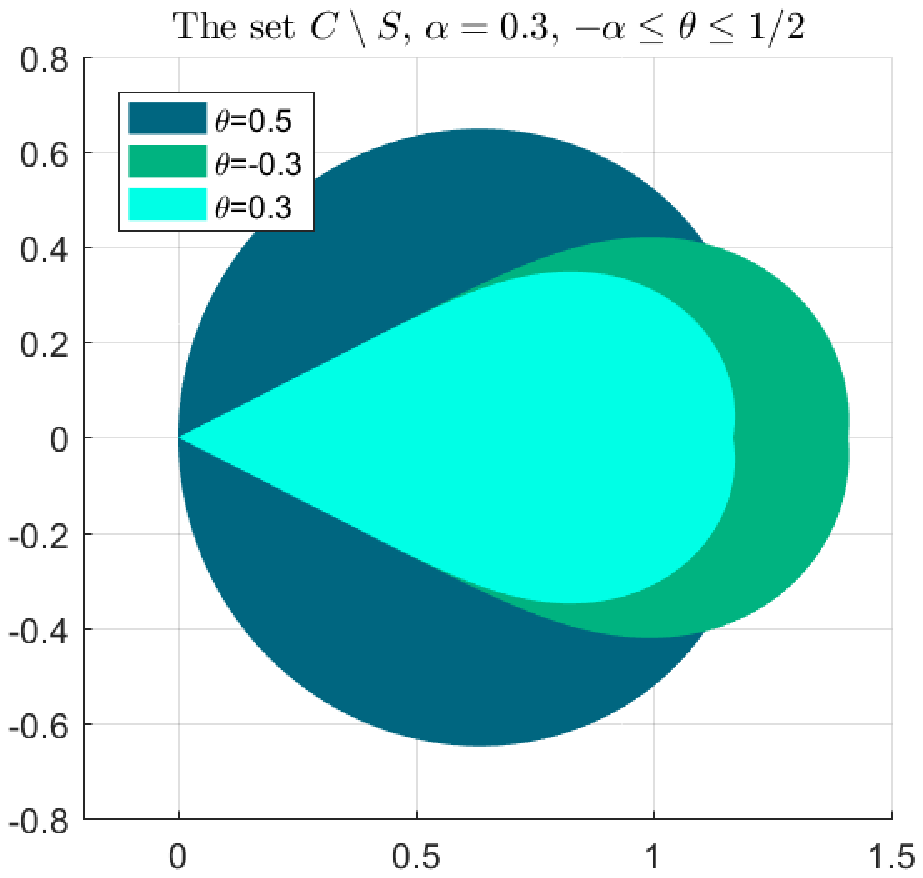}\label{C8}
\end{minipage}
\\
\begin{minipage}{7.5cm}
  \centering\includegraphics[width=7cm]{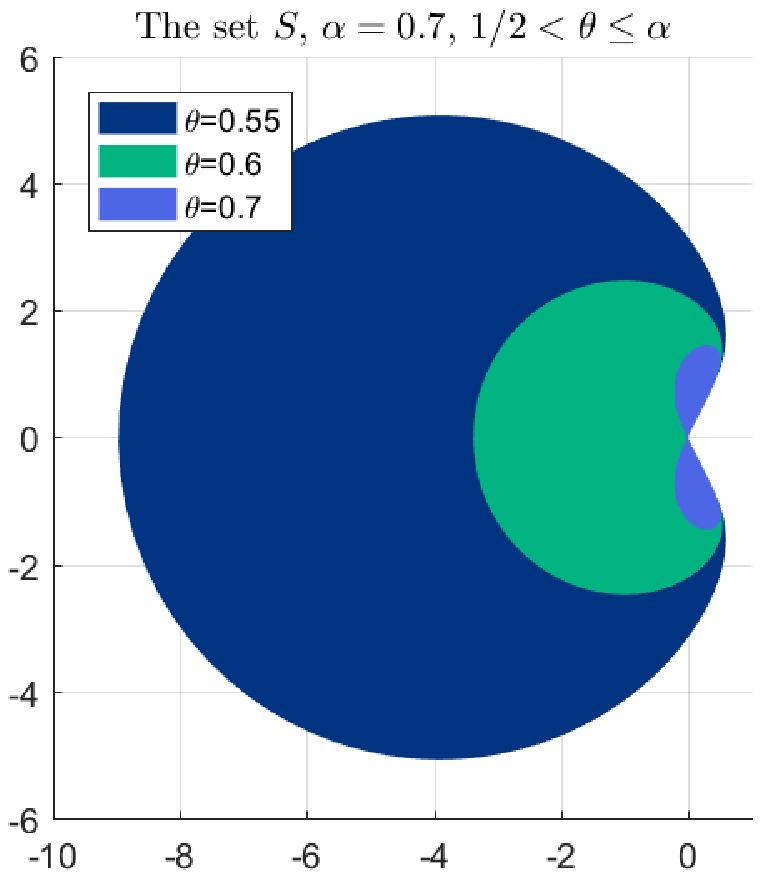}\label{C9}
\end{minipage}
&\begin{minipage}{7.5cm}
  \centering\includegraphics[width=7cm]{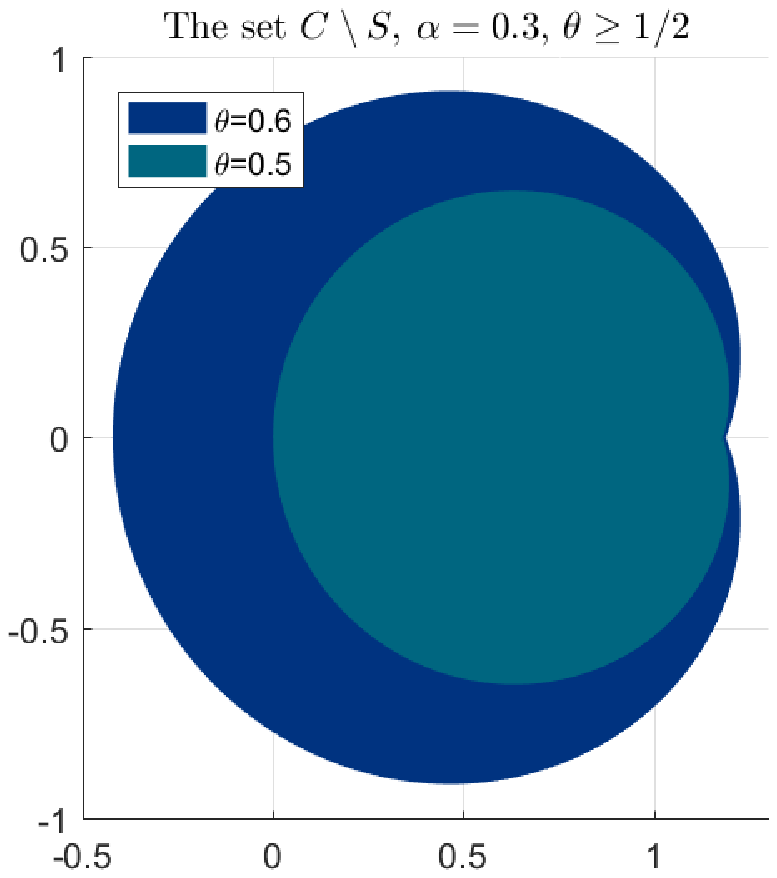}\label{C10}
\end{minipage}
\\
\hline
\end{tabular}}
\end{table}
\begin{table}[h]
\centering
 \caption{The fractional central difference method (\ref{S.41})}\label{tab4}
{\renewcommand{\arraystretch}{2}
\renewcommand{\tabcolsep}{0cm}
\begin{tabular}{|c|c|}
\hline
Scheme-I & Scheme-II \\
\hline
$\omega(\xi)=2^{-\alpha}(1-\xi^2)^{\alpha}$, with $\theta=\alpha$ & $\omega(\xi)=2^{\alpha}(1-\xi^2)^{-\alpha}$, with $\theta=-\alpha$ \\
\hline
\begin{minipage}{7.5cm}
  \centering\includegraphics[width=7cm]{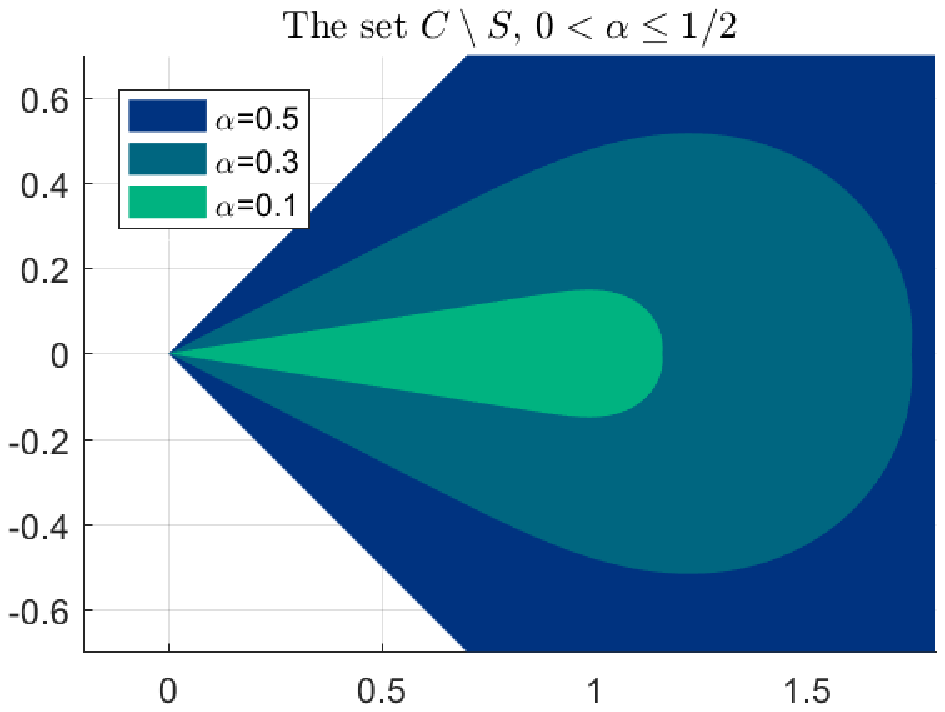}\label{C11}
\end{minipage}
& \begin{minipage}{7.5cm}
  \centering\includegraphics[width=7cm]{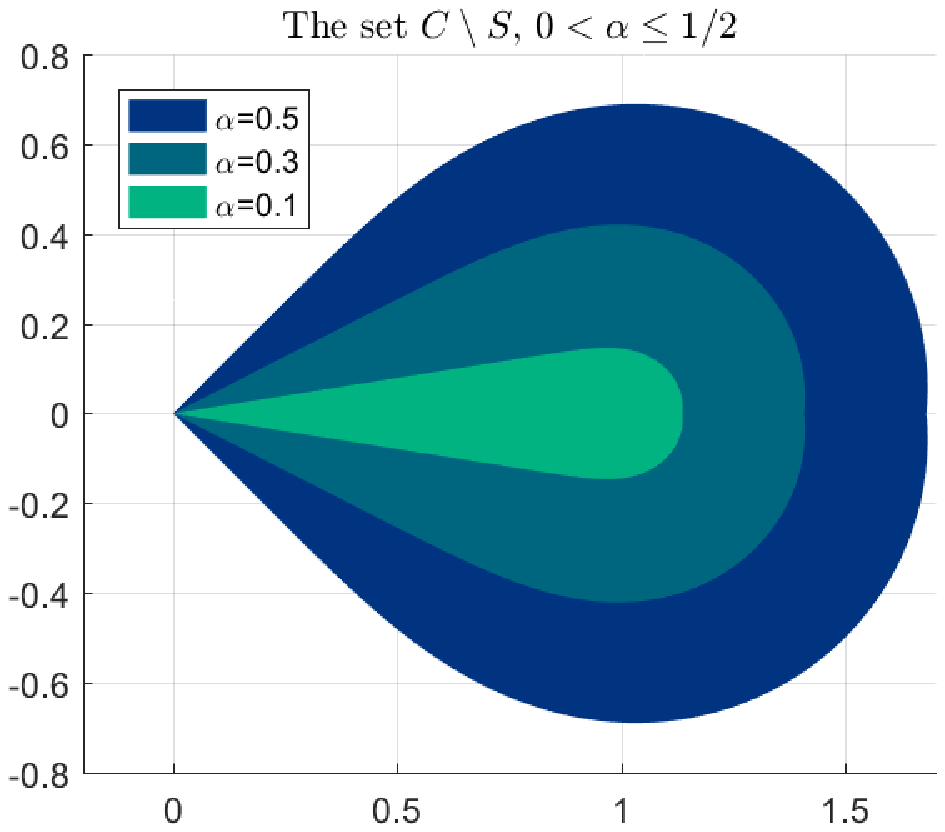}\label{C12}
\end{minipage} \\
\begin{minipage}{7.5cm}
  \centering\includegraphics[width=7cm]{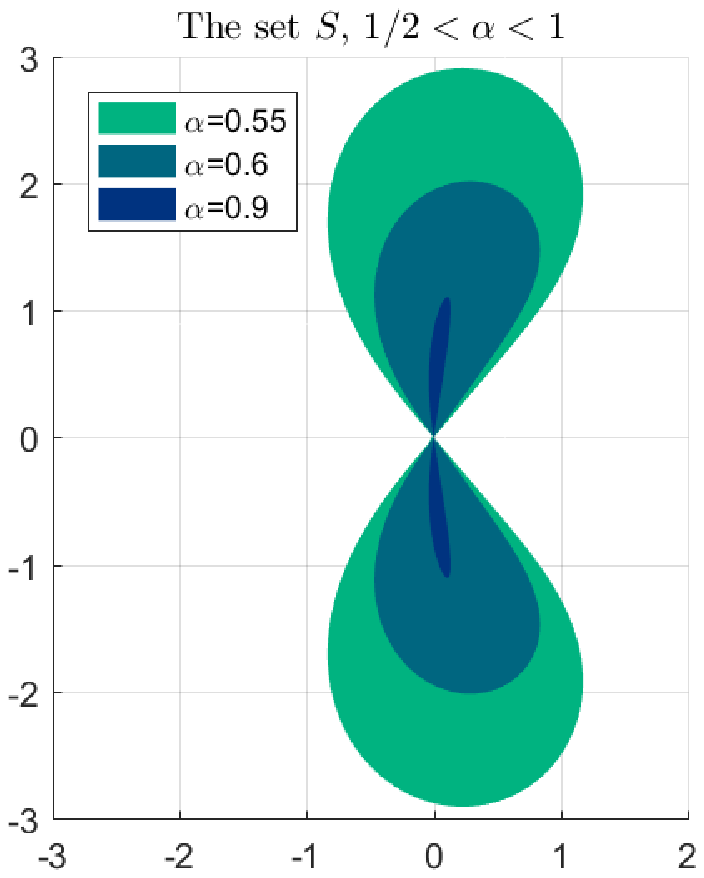}\label{C13}
\end{minipage}
&
\begin{minipage}{7.5cm}
  \centering\includegraphics[width=7cm]{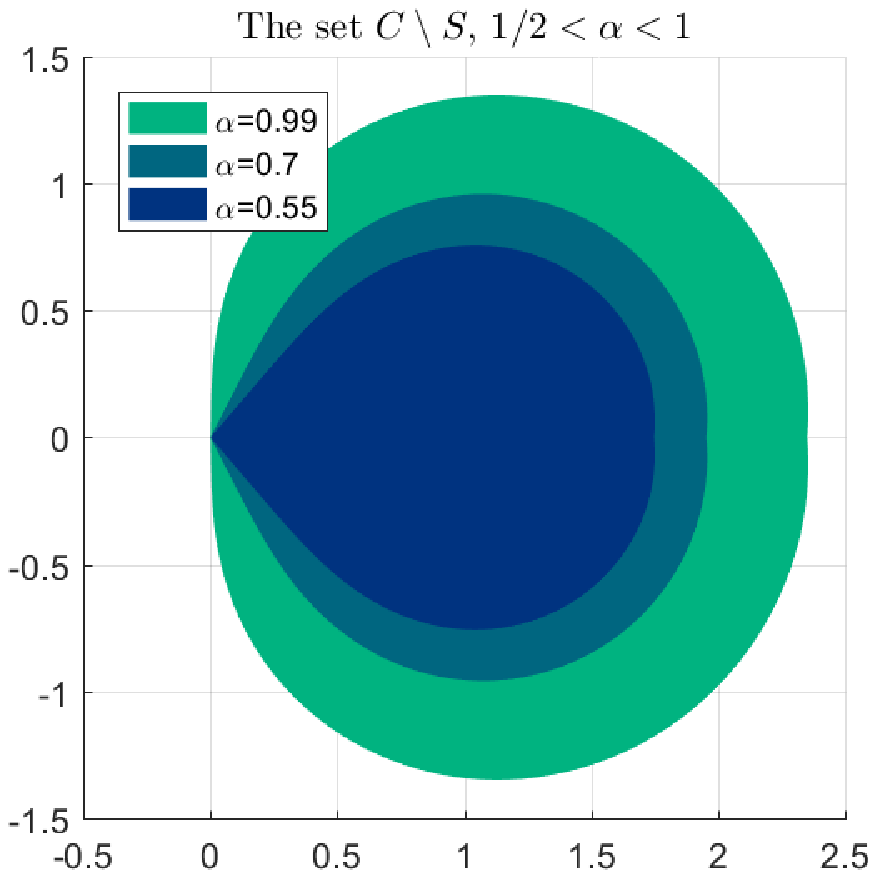}\label{C14}
\end{minipage}\\
\hline
\end{tabular}}
\end{table}
\section{Applications to PDEs}
In this section we apply a class of novel numerical schemes to the time-fractional diffusion equation (\ref{A.1}).
To facilitate the numerical analysis below, we require the scheme satisfies the assumptions of the discrete fractional Gr$\ddot{\text{o}}$nwall inequality \cite{Liao1}.
In the spacial direction, the finite element method is adopted to formulate the fully discrete scheme.
The equation is
\begin{equation}\label{A.1}
\begin{split}
{}_C D_{0,t}^{\alpha} u(t)&=\Delta u+g(x,t)\quad\text{ for } (x,t) \in \Omega \times (0,T] \\
u(x,0)&=u_0(x)~~~~\quad\quad\quad\text{for } x \in \Omega\\
u(x,t)&=0~\quad\quad\quad\quad\quad\quad\text{for } (x,t) \in \partial\Omega \times (0,T),
\end{split}
\end{equation}
where $T>0$, ${}_C D_{0,t}^{\alpha} (\alpha \in (0,1))$ is the Caputo fractional derivative operator defined by (\ref{R.3}), and $\Omega$ is a bounded interval.
The temporal interval is separated uniformly with $0=t_0<t_1<\cdots<t_N=T$ and define $\tau=T/N$ as the step size of the time mesh.
For $\Omega$, we denote $\mathcal{T}_h$ as its shape-regular and quasi-uniform triangulation with the mesh size $h$.
Introduce the subspace $X_h$ of $H_0^1(\Omega)$ defined by
\begin{equation}\label{A.2}\begin{split}
X_h=\{\chi \in H^{1}_0(\Omega):\chi|_{e} \in \mathbb{P}_r(x), e \in \mathcal{T}_h\},
\end{split}
\end{equation}
where $\mathbb{P}_r(x)$ is defined as the set of polynomials (of $x$) with the degree at most $r$$(r \in \mathbb{Z}^{+})$.
\subsection{Numerical scheme}
We approximate the equation (\ref{A.1}) at $t_{n-\theta}$ by the shift-generalized Newton-Gregory formula of order 2 (which is defined by (\ref{S.31}) with $\alpha$ replaced by $-\alpha$), with the generating function
\begin{equation}\label{A.3}
\begin{split}
\omega(\xi)=(1-\xi)^{\alpha}\bigg[1+\bigg(\frac{\alpha}{2}-\theta\bigg)(1-\xi)\bigg] \quad \text{for } \max\big\{0,\frac{\alpha}{2}-\frac{1-\alpha}{1+\alpha}\big\} \leq \theta \leq \frac{\alpha}{2}.
\end{split}
\end{equation}
Considering the relation (\ref{R.4}), we define $v(t):=u(t)-u_0$.
By Theorem \ref{th.2}, we can easily formulate the temporal semi-discrete scheme for (\ref{A.1}) as follows
\begin{equation}\label{A.4}
\begin{split}
\tau^{-\alpha}\sum_{k=0}^{n}\omega_{n-k}v^k=\Delta v^{n-\theta}+g^{n-\theta}+\Delta u_0 +S_h^0\Delta v-S_h^{-\alpha}v,
\end{split}
\end{equation}
where $\Delta v^{n-\theta}:=(1-\theta)\Delta v^n+\theta \Delta v^{n-1}$ and
$g^{n-\theta}:=g(t_{n-\theta})$.
\\
With the space $X_h$, then the fully discrete scheme of (\ref{A.1}) is to find $V^n: [0,T]\to X_h$ such that
\begin{equation}\label{A.5}
\begin{split}
\tau^{-\alpha}\sum_{k=0}^{n}\omega_{n-k}(V^k,\chi_h)+(\nabla V^{n-\theta},\nabla \chi_h)
=(g^{n-\theta}+\Delta u_{0},\chi_h)
-S_h^0(\nabla V,\nabla \chi_h)
-S_h^{-\alpha}(V,\chi_h)
\end{split}
\end{equation}
holds for any $\chi_h \in X_h$.
\par
If we omit the starting part in (\ref{A.5}), we get the scheme
\begin{equation}\label{A.6}
\begin{split}
\tau^{-\alpha}\sum_{k=0}^{n}\omega_{n-k}(V^k,\chi_h)+(\nabla V^{n-\theta},\nabla \chi_h)
=(g^{n-\theta}+\Delta u_{0},\chi_h).
\end{split}
\end{equation}
\subsection{Stability and convergence analysis}
To analyse the stability and convergence of the scheme (\ref{A.1}), we employ the key tool of the discrete fractional Gr$\ddot{\text{o}}$nwall inequality \cite{Liao1}.
In the following analysis, we omit the starting part from our scheme as results will not be affected under suitable conditions.
To meet the assumptions of the discrete fractional Gr\"{o}nwall inequality, we first in Lemmas \ref{l.4}-\ref{l.5} prove some properties of the coefficients of $\omega(\xi)$ defined by (\ref{A.3}).
\par
For clarity, we denote by $\|\cdot\|$ the norm of the $L^2(\Omega)$ space, and by $\|\cdot\|_{s}$ the norm of the Sobolev space $H^s(\Omega)$.
\begin{lemma}\label{l.4}
The coefficients $\omega_k$ defined by the generating function $\omega(\xi)$ in (\ref{A.3}) satisfy
\par
(i) $\omega_0=1+\frac{\alpha}{2}-\theta>0$,\quad $\omega_k\leq 0,~ \text{for } k=1,2,\cdots$, if $\theta \in \big[\frac{\alpha}{2}-\frac{1-\alpha}{1+\alpha},\frac{\alpha}{2}\big]$
\par
(ii) $\sum_{k=0}^{n}\omega_k >0$, and $\sum_{k=0}^{n}\omega_k=\frac{n^{-\alpha}}{\Gamma(1-\alpha)}+O(n^{-1-\alpha})$.
\end{lemma}
\textbf{Proof. } Denote by $\kappa_k$ the coefficients of $(1-\xi)^{\alpha}$.
Then the coefficients $\kappa_k$ satisfy (see \cite{Lubich1,Zeng2}),
\par
(a) $\kappa_n=O(n^{-\alpha-1})$,
\par
(b) $\kappa_0=1, \kappa_n < 0, |\kappa_{n+1}|<|\kappa_n|, n=1,2,\cdots$,
\par
(c) $\displaystyle\sum_{k=0}^{n-1}\kappa_k=\frac{\Gamma(n-\alpha)}{\Gamma(1-\alpha)\Gamma(n)}=\frac{n^{-\alpha}}{\Gamma(1-\alpha)}+O(n^{-1-\alpha}), n=1,2,\cdots.$
\\
On the other hand, we can express $\omega_k$ by $\kappa_k$ as
\begin{equation}\label{A.7}
\begin{split}
\omega_k=\bigg(1+\frac{\alpha}{2}-\theta\bigg)\kappa_k+\bigg(\theta-\frac{\alpha}{2}\bigg)\kappa_{k-1},~ k=1,2,\cdots,
\end{split}
\end{equation}
and, $\omega_0=1+\frac{\alpha}{2}-\theta$.
If $k=1$, we have
\begin{equation}\label{A.8}
\begin{split}
\omega_1=\theta(\alpha+1)-\frac{\alpha}{2}(\alpha+3) \leq -\alpha<0.
\end{split}
\end{equation}
When $k \geq 2$, by (\ref{A.7}) and (b) the condition $\omega_k \leq 0$ is equivalent to the following inequality
\begin{equation}\label{A.8.1}
\begin{split}
\frac{\alpha}{2}-\theta \leq \frac{\kappa_k}{\kappa_{k-1}-\kappa_k}.
\end{split}
\end{equation}
With the estimates
\begin{equation}\label{A.8.2}
\begin{split}
\frac{\kappa_{k-1}}{\kappa_k}=\frac{k}{k-\alpha-1}=\frac{1}{1-\frac{\alpha+1}{k}}\leq \frac{1}{1-\frac{\alpha+1}{2}},
\end{split}
\end{equation}
we can easily check that (\ref{A.8.2}) holds for $\theta \in \big[\frac{\alpha}{2}-\frac{1-\alpha}{1+\alpha},\frac{\alpha}{2}\big]$.
\\
By the definition of $\omega(\xi)$, we know that $\sum_{k=0}^{\infty}\omega_k=0$, hence by (i), we can get $\sum_{k=0}^{n}\omega_k > 0$.
Now combining (a), (c) with (\ref{A.7}), we have
\begin{equation}\label{A.9}
\begin{split}
\sum_{k=0}^{n-1}\omega_k=\sum_{k=0}^{n-1}\kappa_k+\bigg(\frac{\alpha}{2}-\theta\bigg)\kappa_{n-1}
=\frac{n^{-\alpha}}{\Gamma(1-\alpha)}+O(n^{-1-\alpha}).
\end{split}
\end{equation}
The proof of the lemma is completed.
\begin{lemma}\label{l.5}
Let $A_{k}:=\frac{1}{\tau^{\alpha}}\sum_{j=0}^{k}\omega_j$.
Then, $A_k \geq A_{k+1}>0 (k \geq 0)$.
Furthermore, there exists a positive $\pi_A$ such that
\begin{equation}\label{A.10}
\begin{split}
A_{n-k}\geq \frac{1}{\pi_A \tau}\int_{t_{k-1}}^{t_k} \frac{(t_n-s)^{-\alpha}}{\Gamma(1-\alpha)}\mathrm{d}s,  \quad \text{ for any } n \geq k \geq 1.
\end{split}
\end{equation}
\end{lemma}
\textbf{Proof. } By careful derivation one can see that to prove the existence of $\pi_A$, it is sufficient to demonstrate
\begin{equation}\label{A.11}
\begin{split}
\pi_A \geq \frac{1}{\Gamma(2-\alpha)}\max_{n \geq 0} \Theta_n >0,
\end{split}
\end{equation}
where $\Theta_n=\frac{(n+1)^{1-\alpha}-n^{1-\alpha}}{\sum_{j=0}^{n}\omega_j}$.
By Lemma \ref{l.4} (ii), $\Theta_n$ is an increasing sequence, and the limit $\lim_{n\to \infty}\Theta_n$ exists.
The proof is completed.
\par
The complementary discrete convolution kernels $P_j$ is essential to the development of the inequality (\ref{A.15}), which is defined as the coefficients of the function $P(\xi)=\tau^{\alpha}\omega(\xi)^{-1}$.
We remark that this definition implies
\begin{equation}\label{A.12}
\begin{split}
\sum_{j=0}^{n}P_{n-j}A_j \equiv 1 \quad \text{for any } n \geq 0.
\end{split}
\end{equation}
Actually, the sequence $\{A_k\}$ defined in Lemma \ref{l.5} is the coefficients of the function $A(\xi)=\frac{\omega(\xi)}{\tau^{\alpha}(1-\xi)}$.
On the other hand,
$\sum_{j=0}^{n}P_{n-j}A_j$ is the $n$th coefficient of $P(\xi)A(\xi)=\frac{1}{1-\xi}=1+\xi+\xi^2+\cdots$, which means (\ref{A.12}).
\par
Now based on Lemmas \ref{l.4}-\ref{l.5}, for coefficients $A_k$ defined in Lemma \ref{l.5}, we have the following discrete fractional Gr$\ddot{\text{o}}$nwall inequality, see \cite{Liao1}.
\begin{lemma}\label{l.6}
Let $\{\phi^n\}_{n=1}^{\infty}$ and $\{\eta_l\}_{l=0}^{\infty}$ be given nonnegative sequences.
Let $0 \leq \theta <1$.
Assume further the series $\sum_{l=0}^{\infty}\eta_l$ is bounded with $\Lambda$, i.e., $\sum_{l=0}^{\infty}\eta_l < \Lambda$, and that the time step size satisfies
\begin{equation}\label{A.13}
\begin{split}
\tau \leq \frac{1}{{}^{\alpha}\sqrt{2\pi_A \Gamma(2-\alpha)\Lambda}}.
\end{split}
\end{equation}
Then for any nonnegative sequence $\{v^k\}_{k=0}^N$ such that
\begin{equation}\label{A.14}
\begin{split}
\sum_{k=1}^{n}A_{n-k}\big[(v^k)^2-(v^{k-1})^2\big]
\leq \sum_{k=1}^{n}\eta_{n-k}(v^{k-\theta})^2+v^{n-\theta}\phi^n \quad \text{for } 1 \leq n \leq N,
\end{split}
\end{equation}
it holds that
\begin{equation}\label{A.15}
\begin{split}
v^n \leq 2 E_{\alpha}(2\pi_A \Lambda t_n^{\alpha})\bigg(v^0+\sum_{j=1}^{n}P_{n-j}\phi^j\bigg) \quad \text{for } 1 \leq n \leq N.
\end{split}
\end{equation}
\end{lemma}
\begin{rem}\label{rem.4}
A careful examination shows that $\sum_{j=1}^{n}P_{n-j}$ is bounded by some constant $C$ which is independent of $n$.
Hence, (\ref{A.15}) can be simplified to
\begin{equation}\label{A.16}
\begin{split}
v^n \leq C(v^0+\max_{1\leq j \leq n}\phi^j) \quad \text{for } 1 \leq n \leq N.
\end{split}
\end{equation}
\end{rem}
We reformulate the numerical scheme (\ref{A.6}) with the coefficients $A_k$ defined in Lemma \ref{l.5} as the following
\begin{equation}\label{A.17}
\begin{split}
\sum_{k=1}^{n}A_{n-k}(V^k-V^{k-1},\chi_h)+(\nabla V^{n-\theta},\nabla \chi_h)
=&(g^{n-\theta}+\Delta u_{0},\chi_h).
\end{split}
\end{equation}
Now we are in a position to analyse the stability and error estimates of the scheme by the same approaches taken in \cite{Liao1}, hence, we omit the proof of the following theorems.
\begin{thm}\label{th.5}
Define $U^n:=V^n+u_{0h}$ and $V^n$ is the numerical solution of (\ref{A.6}), $u_{0h}$ is a proper approximation of $u_0$.
Then with $\tau$ satisfying (\ref{A.13}), we have the stability estimates as
\begin{equation}\label{A.18}
\begin{split}
\|U^n\|\leq \|u_{0h}\|+ C\|\Delta u_{0}\| + C \max_{1\leq j\leq n}\|g^{n}\|,
\end{split}
\end{equation}
where the constant $C$ is independent of $\tau$ and $h$.
\end{thm}
\begin{thm}\label{th.6}
Let $V^n$ be the numerical solution of (\ref{A.6}).
Define $U^n:=V^n+u_{0h}$, where $u_{0h}$ is a proper approximation of $u_0$.
Then $U^n$ is a numerical solution of (\ref{A.1}).
Suppose $u\in C([0,T]; H^{r+1}(\Omega) \cap H_0^1(\Omega)) \cap C^3([0,T]; L^2(\Omega))$, then with $\tau$ satisfying (\ref{A.13}),  we have the error estimates as
\begin{equation}\label{A.19}
\begin{split}
\|U^n-u^n\|\leq C(\tau^2+h^{r+1}),
\end{split}
\end{equation}
where the constant $C$ is independent of $\tau$ and $h$.
\end{thm}
\subsection{Numerical experiments}
In this subsection, we implement some numerical experiments to further confirm our convergence estimates (\ref{A.19}).
To this end, define the error $E(\tau,h):=\max_{0 \leq n \leq N}\|u^n-U^n\|$.
The convergence rate are derived by the formula
\begin{equation}\label{A.20}\begin{split}
\text{temporal order }=\log_2\frac{E(2\tau,h)}{E(\tau,h)},
\quad
\text{spatial order }=\log_2\frac{E(\tau,2h)}{E(\tau,h)}.
\end{split}\end{equation}
Let $\Omega=(0,1)$, $T=1$ and divide the interval $\Omega$ as $0=x_0<x_1<\cdots<x_M=1$ with the mesh size $h=1/M$.
Define the finite element space $X_h$ as the set of piecewise linear polynomials with $r=1$.
The exact solution is taken as
\begin{equation}\label{A.21}
\begin{split}
u(x,t)=(1+t^{\alpha}+t^{2\alpha}+t^3)\sin(2\pi x),
\end{split}
\end{equation}
and the corresponding source term $g$ is
\begin{equation}\label{A.22}
\begin{split}
g(x,t)=\sin(2\pi x)\bigg[\Gamma(1+\alpha)+4\pi^2(1+t^{\alpha}+t^{2\alpha}+t^3)+\frac{2\alpha t^{\alpha}\Gamma(2\alpha)}{\Gamma(\alpha+1)}+\frac{6t^{3-\alpha}}{\Gamma(4-\alpha)}\bigg].
\end{split}
\end{equation}
One can see that there is some weak singularity for $u$ at initial value, hence, by Theorem \ref{th.2} we add the starting part to obtain a second-order convergence in time.
In the following tables, we denoted by $E_c(\tau,h)$ the errors of the scheme with starting part, and by $E_o(\tau,h)$ the errors without the starting part.
\par
In Table \ref{tab5}, with fixed fine space mesh size $h=\frac{1}{5000}$, we take different $\alpha$ and for each $\alpha$ we choose some $\theta$ which satisfy $\max\big\{0,\frac{\alpha}{2}-\frac{1-\alpha}{1+\alpha}\big\} \leq \theta \leq \frac{\alpha}{2}$.
Now one can easily check that with the starting part, we have obtained a second-order convergence in time.
For the scheme without the starting part, the convergence rate is much lower, especially when $\alpha$ is small.
\par
In Table \ref{tab6}, we fix the time mesh size $\tau=\frac{1}{1000}$, and similarly, compare the spacial convergence rate for the scheme with or without the starting part.
We can see that the convergence rate is optimal when the starting part is added, and for those $\alpha$ close to zero (which means there is a stronger singularity for the solution), the rate becomes much lower if the starting part is omitted.
\begin{table}[h]
\centering
 \caption{The temporal convergence rate with $h=\frac{1}{5000}$}\label{tab5}
{\renewcommand{\tabcolsep}{0.3cm}
 \begin{tabular}{|c|c|c|c|c|c|c|}
\hline
  $\alpha$ &$\theta$&  $\tau$ & $E_c(\tau, h)$ & rate &  $E_o(\tau, h)$ & rate \\
\hline
	&		&	  1/10  	&	2.35436E-05	&	---	&	2.87659E-03	&	---	\\
	&	0	&	  1/20  	&	5.47471E-06	&	2.1045 	&	2.79338E-03	&	0.0423 	\\
	&		&	  1/40  	&	1.13371E-06	&	2.2717 	&	2.71399E-03	&	0.0416 	\\
0.1	&		&	  1/80  	&	2.55465E-07	&	2.1499 	&	2.63920E-03	&	0.0403 	\\
\cline{2-7}
	&		&	  1/10  	&	9.59982E-04	&	---	&	4.22249E-02	&	---	\\
	&	0.05	&	  1/20  	&	2.43373E-04	&	1.9798 	&	3.81734E-02	&	0.1455 	\\
	&		&	  1/40  	&	6.15520E-05	&	1.9833 	&	3.44973E-02	&	0.1461 	\\
	&		&	  1/80  	&	1.57177E-05	&	1.9694 	&	3.11936E-02	&	0.1452 	\\
\hline
	&		&	  1/10  	&	1.72824E-03	&	---	&	5.64791E-03	&	---	\\
	&	0.1	&	  1/20  	&	4.38829E-04	&	1.9776 	&	2.56403E-03	&	1.1393 	\\
	&		&	  1/40  	&	1.10750E-04	&	1.9864 	&	4.53468E-04	&	2.4993 	\\
0.5	&		&	  1/80  	&	2.80539E-05	&	1.9810 	&	8.92050E-04	&	-0.9761 	\\
\cline{2-7}
	&		&	  1/10  	&	3.22333E-03	&	---	&	2.00719E-02	&	---	\\
	&	0.2	&	  1/20  	&	8.15120E-04	&	1.9835 	&	1.28495E-02	&	0.6435 	\\
	&		&	  1/40  	&	2.05120E-04	&	1.9905 	&	7.75159E-03	&	0.7291 	\\
	&		&	  1/80  	&	5.16829E-05	&	1.9887 	&	4.28721E-03	&	0.8545 	\\
\cline{2-7}
\hline
	&		&	  1/10  	&	4.70189E-03	&	---	&	6.08467E-03	&	---	\\
	&	0.4	&	  1/20  	&	1.21112E-03	&	1.9569 	&	1.52565E-03	&	1.9958 	\\
	&		&	  1/40  	&	3.07190E-04	&	1.9791 	&	6.99968E-04	&	1.1241 	\\
0.9	&		&	  1/80  	&	7.75617E-05	&	1.9857 	&	6.01301E-04	&	0.2192 	\\
\cline{2-7}
	&		&	  1/10  	&	4.90787E-03	&	---	&	6.31974E-03	&	---	\\
	&	0.45	&	  1/20  	&	1.25984E-03	&	1.9619 	&	1.58234E-03	&	1.9978 	\\
	&		&	  1/40  	&	3.19039E-04	&	1.9814 	&	3.96145E-04	&	1.9980 	\\
	&		&	  1/80  	&	8.04831E-05	&	1.9870 	&	1.45886E-04	&	1.4412 	\\
\hline
 \end{tabular}}
\end{table}
\begin{table}[h]
\centering
 \caption{The spacial convergence rate with $\tau=\frac{1}{1000}$}\label{tab6}
{\renewcommand{\tabcolsep}{0.3cm}
 \begin{tabular}{|c|c|c|c|c|c|c|}
\hline
  $\alpha$ &$\theta$&  $h$ & $E_c(\tau, h)$ & rate &  $E_o(\tau, h)$ & rate \\
\hline
	&		&	  1/10  	&	9.90558E-02	&	---	&	9.90559E-02	&	---	\\
	&	0	&	  1/20  	&	2.49005E-02	&	1.9921 	&	2.49007E-02	&	1.9921 	\\
	&		&	  1/40  	&	6.23369E-03	&	1.9980 	&	6.23387E-03	&	1.9980 	\\
0.3	&		&	  1/80  	&	1.55895E-03	&	1.9995 	&	4.29525E-03	&	0.5374 	\\
\cline{2-7}
	&		&	  1/10  	&	9.90560E-02	&	---	&	9.90561E-02	&	---	\\
	&	0.15	&	  1/20  	&	2.49008E-02	&	1.9921 	&	2.49009E-02	&	1.9920 	\\
	&		&	  1/40  	&	6.23393E-03	&	1.9980 	&	6.23403E-03	&	1.9980 	\\
	&		&	  1/80  	&	1.55920E-03	&	1.9993 	&	7.03535E-03	&	-0.1745 	\\
\hline
	&		&	  1/10  	&	9.85495E-02	&	---	&	9.85515E-02	&	---	\\
	&	0.3	&	  1/20  	&	2.47681E-02	&	1.9924 	&	2.47686E-02	&	1.9924 	\\
	&		&	  1/40  	&	6.20048E-03	&	1.9980 	&	6.20066E-03	&	1.9980 	\\
0.8	&		&	  1/80  	&	1.55092E-03	&	1.9993 	&	1.55102E-03	&	1.9992 	\\
\cline{2-7}
	&		&	  1/10  	&	9.85493E-02	&	---	&	9.85515E-02	&	---	\\
	&	0.4	&	  1/20  	&	2.47681E-02	&	1.9924 	&	2.47687E-02	&	1.9924 	\\
	&		&	  1/40  	&	6.20052E-03	&	1.9980 	&	6.20072E-03	&	1.9980 	\\
	&		&	  1/80  	&	1.55098E-03	&	1.9992 	&	1.55108E-03	&	1.9992 	\\
\hline
 \end{tabular}}
\end{table}
\section{Conclusion}
In this paper, the shifted convolution quadrature theory is developed based on the extensible framework established by Lubich.
The definition of consistency is generalized and the equivalent theorem is established for the $SCQ$ theory.
The fertility of the generalized framework is demonstrated by being able to transform generating functions of the $CQ$ to those of the $SCQ$ (Theorem {\ref{th.3}}), develop the shift-generalized Newton-Gregory formula (Corollary {\ref{cor.2}}), include as many existing approximation formula with integer convergence rate as possible (Sec. 2) and easy to design new formulas with desired structure (Example 4).
We shall point out that the allowable structure of a generating function is not limited to the case $(\ref{S.2.1})$, as Lubich \cite{Lubich1} said "condition (\ref{S.2.1}) can be considerably relaxed, however, the class (\ref{S.2.1}) is probably large enough for all practical applications".
Indeed, all formulas mentioned in this paper have further assumed $p_4(\xi)\equiv1$ in (\ref{S.2.1}), and it is not hard to construct formulas with $p_4(\xi)\neq 1$ identically which are stable and consistent.
Another merit of the $SCQ$ is that it inherits the correction techniques from the $CQ$, by which the high-order convergence rate can be obtained numerically despite the weak regularity for solutions at initial node.
To further explore the stability properties of the schemes proposed, we analyse the impact of the parameter $\theta$ on the stable regions with changing $\theta$.
We emphasize that the shift parameter $\theta$ plays an essential role in the stability of a numerical scheme, and careful examination for the choice of $\theta$ is of vital importance before employing new schemes.
For some special designed generating functions, we apply the resulted schemes to PDEs with the technique of the fractional Gr\"{o}nwall inequality to analyse the stability and convergence.
The results of the numerical experiments further confirm our theory analysis.
\par
Authors think there are at least two approaches in our future work:
a) Propose some simple generating functions that meet with the techniques already developed for the numerical analysis for different PDEs, and
b) develop new techniques that are suitable for as many different types of generating functions as possible.
\section*{Acknowledgements}
The work of the first author was supported in part by the NSFC grant 11661058.
The work of the third author was supported in part by the NSFC grant 11761053, the NSF of Inner Mongolia 2017MS0107, and the program for Young Talents of Science and Technology in Universities of Inner Mongolia Autonomous Region NJYT-17-A07.
The work of the fourth author was supported in part by grants NSFC 11871092 and NSAF U1530401.

\end{document}